\documentclass[a4paper,reqno]{amsart}
\usepackage[english]{babel}
\usepackage{amsfonts, amssymb, mathtools}
\usepackage[a4paper]{geometry}
\usepackage[toc]{appendix}
\usepackage[final]{hyperref}
\newtheorem{theorem}{Theorem}[section]
\newtheorem{lemma}{Lemma}[section]
\newtheorem{corollary}{Corollary}[section]
\newtheorem{proposition}{Proposition}[section]

\theoremstyle{definition}
\newtheorem{definition}{Definition}[section]

\theoremstyle{remark}
\newtheorem{remark}{Remark}[section]
\newtheorem{example}{Example}[section]

\numberwithin{equation}{section}

\allowdisplaybreaks

\newcommand{\e}{\mathrm{e}}
\newcommand{\1}{1\!\!\,{\textrm I}}
\newcommand{\ov}{\overline}

\newcommand{\wt}{\widetilde}
\newcommand{\wh}{\widehat}
\newcommand{\pt}{\partial}
\newcommand{\Id}{\mathrm{Id}}
\newcommand{\Law}{\mathrm{Law}}
\newcommand{\Leb}{\mathrm{Leb}}

\newcommand{\cadlag}{c\`adl\`ag }
\newcommand{\caglad}{c\`agl\`ad }

\DeclareMathOperator{\E}{E}

\DeclareMathOperator{\Prob}{P}
\DeclareMathOperator{\Var}{Var}

\newcommand{\vf}{\varphi}
\newcommand{\vk}{\kappa}
\newcommand{\ve}{\varepsilon}

\newcommand{\mbR}{{\mathbb R}}

\newcommand{\mbN}{{\mathbb N}}

\newcommand{\cF}{{\mathcal F}}

\newcommand{\cN}{{\mathcal N}}

\newcommand{\cA}{{\mathcal A}}
\newcommand{\cM}{{\mathcal M}}

\newcommand{\cW}{{\mathcal W}}
\newcommand{\cR}{{\mathcal R}}

\newcommand{\mcC}{{\mathcal C}}

\makeatletter
\@namedef{subjclassname@2020}{
  \textup{2020} Mathematics Subject Classification}
\makeatother

\title[Splitting for some classes of stochastic flows]{Splitting for some classes of homeomorphic and coalescing stochastic flows}
\author{M.B. Vovchanskyi}
\address{Institute of Mathematics, National Academy of Sciences of Ukraine, Teresh\-chenkivska Str. 3, Kyiv 01601, Ukraine}
\email{vovchansky.m@gmail.com}


\usepackage[backend=biber, 
	bibencoding=utf8,
	style=ieee, 
	citestyle=numeric-comp,
	sorting=nty,
	autolang=other,
	maxbibnames=10,
	giveninits=true,
	date=year,	doi=false,isbn=false,url=false,
	eprint=true,
	hyperref=true
	]{biblatex}
\AtEveryBibitem{\clearfield{language}}
\AtEveryBibitem{\clearfield{note}}

\begin{document}

\begin{abstract}
The splitting scheme (the Kato-Trotter formula) is applied  to stochastic flows with common noise of the type introduced  by Th.E.~Harris. The case of possibly coalescing flows with continuous infinitesimal covariance is considered and the weak convergence of the corresponding finite-dimensional motions is established. As applications, results for the convergence of the associated pushforward measures and dual flows are given. Similarities between splitting and the Euler-Maruyama scheme yield estimates of the speed of the convergence under additional regularity assumptions.
\end{abstract}

\maketitle 

\section{Introduction}

Introduced (with zero drift) in \cite{Ha84Coalescing}, Harris flows are families of random transformation of the real line that represent the joint movement of one-dimensional interacting Brownian particles whose pairwise correlation depends on the distance between them and is given via so-called infinitesimal covariance $\vf.$   
Since coalescence is possible, such transformations are not, in contrast to the case of diffeomorphic flows obtained as solutions to SDEs with sufficiently smooth coefficients~\cite{Ku90Stochastic}, necessarily continuous. One natural and straightforward extension of the notion of the Harris flow is to add drift to affect the motion of particles in a way similar to the case of the Arratia flow in \cite{Do07MeasureEng}. This brings us closer to biological and physical models that use potentials of different forms \cite{Bress19Stochastic, KoLeitMann09Correlated}  while introducing common noise as in \cite{CoFlan16Propagation, GuLu21ScalingArxiv} including one that forces particles to collide. 

The main goal of the paper is to apply the well-known  method of splitting in the stochastic setting~\cite{GyKry03Splitting, Faou09Analysis, GyRa11Note, GonKo98Fractional, BenGloRas92Approximations, BuckSamTam22Splitting, BreGou19Analysis, CuiHong19strong, BreCuiHong19strong} to Harris flows, so that the actions of the semigroups generated by the corresponding driftless Harris flow and the ordinary ODE are separated.

The formal definition of a Harris flow with drift adopted in the paper is based of the definition of a driftless Harris flow from \cite{WarWa04Spectra} (see also \cite{Ha84Coalescing, Mat89Coalescing, AmaTaYu19Convergence}), with a minor modification as in \cite{Vov18Convergence}. Let $D^{\uparrow}(\mbR)$ be a separable metrizable topological space of non-decreasing \cadlag functions on $\mbR$ equipped with the Skorokhod $J_1$ topology~\cite{Bi68Convergence, EthKur09Markov}. 
Since for any $f, g\in D^\uparrow(\mbR)$ the composition $f\circ g \in D^{\uparrow}(\mbR)$~\cite[Lemma 13.2.4]{Whitt02Stochastic}, and, for $ D^\uparrow(\mbR)-$valued random elements $\xi, \eta,$ we have
\[
\left\{ \xi \circ \eta (t) \ge a \right\} = \left\{ \eta(t) \ge \hat\xi(a) \right\},
\]
where $\hat f$ is the \caglad generalized inverse  of $f,$ 
the composition $\xi \circ \eta$  defines a random element in $D^\uparrow(\mbR).$ The space of $\mbR^d-$valued \cadlag functions with non-decreasing coordinates  endowed with the $J_1-$topology is denoted by $D^\uparrow(\mbR^d),$ and  standard Skorokhod spaces of functions on $[0,T], T >0,$ with values in $\mbR^d$ are denoted $D([0,T],\mbR^d).$ 

\begin{remark}
Standard sources discuss spaces $D([0,\infty), \mbR^d)$ or $D([0,T), \mbR^d)$ but the extension of the parametric set to the whole $\mbR$ can be found in~\cite{FerVo09Weak}.   
\end{remark}

\begin{definition}
\label{harris.flow}
A Harris flow $X^{\vf,a}$ with the infinitesimal covariance $\vf$ and drift $a$ is a family of $D^\uparrow(\mbR)-$valued random elements $\{ X^{\vf,a}_{s,t}(\cdot) \mid 0\le s\le t\}$ such that
\begin{enumerate}
\item
for any $s\le t\le r,$ $\Prob\left( X^{\vf,a}_{s,r}=X^{\vf,a}_{t,r}\circ X^{\vf,a}_{s,t}\right)=1$; $X^{\vf,a}_{s,s}=\Id$  a.s. ($\Id$ is the identity mapping);
\item
for any $t_1\le t_2 \le \ldots \le t_n$ random elements $X^{\vf,a}_{t_1,t_2},\ldots,X^{\vf,a}_{t_{n-1},t_n}$ are independent;
\item
for any $s, t, h>0, s< t,$ $\Law\left(X^{\vf,a}_{s,t}\right)=\Law\left(X^{\vf,a}_{s+h,t+h}\right);$
\item
as $h\rightarrow 0+,$ $X^{\vf,a}_{0,h}\rightarrow \Id$ in probability in $D^\uparrow(\mbR);$
\item
for any $x \in\mbR, s\ge 0,$ the  process 
\[
t\mapsto w_t(x,s) \equiv X^{\vf,a}_{s,t}(x) - x - \int_s^t a(X^{\vf,a}_{s,r}(x)) dr, \quad t \ge s, 
\]
is a $(\cF^{X^{\vf,a}}_{s,r})_{r\ge s}-$Wiener process started at $0,$ where 
\[
\cF^{X^{\vf,a}}_{s,t} = \sigma\left\{ X^{\vf,a}_{u_1,u_2}, s \le  u_1\le u_2 \le t\right\}, \quad 0 \le s \le t;
\]  
\item
for any $x,y\in\mbR, s \ge 0$ 
\[
\left\langle w(x,s), w(y,s) \right\rangle_t = \int_s^t \vf\left(X^{\vf,a}_{s,r}(x)-X^{\vf,a}_{s,r}(y)\right)dr, \quad t \ge s.
\]
\end{enumerate} 
\end{definition}

The following splitting scheme is used. Let $0 =t_0 < t_1 <\ldots < t_m = T$  for some $T.$ Define recursively  piecewise continuous processes $(u,y)$ such that for any $k=\ov{0,m-1}$ 
\begin{align*}
u_t &=  y_{t_k-} + \int_{t_k}^t a\left( u_s \right) ds, \\
y_t &= X^{\vf,a}_{t_k, t} \left( u_{t_{k+1}-}\right), \\
& \qquad t \in [t_k, t_{k+1}),
\end{align*} 
with $y_{0-}=x.$

We establish the weak convergence of finite-dimensional motions in Skorokhod spaces in the general case of continuous $\vf$ as the size of a partition tends to $0$ (Theorem~\ref{th:con.1}). This result is used to derive the convergence of the pushforward measures under the actions of the corresponding flows under an additional assumption that guarantees the initial flow to be a coalescing one (Theorem~\ref{th:add.1}). As a second application, the convergence of the associated dual flows in the reversed time is established (Theorem~\ref{th:dual.1}). 

If a Harris flow admits a representation as the unique strong solution of an SDE (see Section~\ref{sec:sdes} for details), which corresponds to the additional  assumption of $\sqrt{1-\vf}$ being Holder continuous of order $\beta\ge \frac{1}{2},$ one can (almost) mechanically transfer proofs and conclusions for the Euler-Maruyama scheme~\cite{BenGloRas92Approximations, GyRa11Note} into our setting (Theorem~\ref{th:homeo.1}). We emphasize the highly derivative nature of the results in this case. To formulate the results, the Wasserstein distance on the space of distributions of random measures is chosen (Theorem~\ref{th:homeo.2}).

\section{Existence of Harris flows with non-trivial drift}
Sufficient conditions for $X^{\vf,a}$ to exist and be a coalescing flow are given in this section.

Hereinafter $\mbR_+ = (0,\infty).$
\begin{definition}
\label{definition.vf}
A continuous symmetric function $\vf\colon \mbR \mapsto \mbR$ belongs to $\Phi_{*}$ if 
\begin{enumerate}
\item $\vf$ is strictly positive definite; 
\item $\vf(0)=1;$
\item $\vf$ is Lipschitz continuous outside any neighborhood of $0.$ 
\end{enumerate}
\end{definition}

\begin{definition}
\label{definition.vf.cond}
A function $\vf\in\Phi_{*}$ belongs to $\Phi_{\alpha}$ for some $\alpha \in (0,2]$ if for some $ C_\vf, \wt{C}_\vf > 0$ 
\begin{equation}
\label{eq:phi}
1-\vf(x) \ge C_\vf |x|^\alpha, \quad x \in [-\wt{C}_\vf,\wt{C}_\vf]. 
\end{equation}
\end{definition}

\begin{definition}
\label{definition.drift}
A measurable function $a\colon \mbR \mapsto \mbR$ belongs to $A_{\beta}$ for some $\beta\in\mbR$ if there exists a non-negative function  $ \rho \in C(\mbR_+)$ such that for some $C_\rho, \wt{C}_\rho > 0$
\begin{align*}
|a(x+y) - a(x)| & \le  \rho(y), \quad x \in \mbR, y\in\mbR_+, \\
\rho(x) &\le C_\rho x^\beta, \quad x \in (0, \wt{C}_\rho], \\
 \rho(x) &\le C_\rho(1 + |x|), \quad x \in \mbR_+.
\end{align*}
\end{definition}

\begin{remark}
If $\vf \in \Phi_*,$ then $\vf(x) < 1$ when $x\not=0$ and $\vf$ is bounded. A function $a\in A_\beta$ does not have to be right continuous or bounded at $0$ if $\beta < 0.$
\end{remark}

\begin{example}
If $\vf(x) = \e^{-|x|^\alpha}, x\in\mbR,$ for $\alpha\in (0,2],$ then $\vf\in \Phi_{\alpha}.$ Here $\alpha=2$ corresponds to a homeomorphic non-coalescing flow. 
\end{example}

\begin{example}
The Brownian web (the Arratia flow) \cite{FonIsoNewRa04Brownian} can be seen as an extreme example of the Harris flow with $\vf(x)=\1[x=0].$ In this case particles are independent before a collision.
\end{example}

\begin{example}
Consider the following example of $\vf\in\Phi_{3/2}$ that does not satisfy the Holder condition of any positive order. Given positive $C_1, C_2$ with $C_1+C_2=1$ set
\begin{align*}
\vf(x) &= C_1\e^{-\frac{x^2}{2}} + C_2 \vf_1(x), \\
\vf_1(x) &= \sum_{n\ge 1} \frac{\cos(\e^{n/2} n^{3/2} x)}{n^2}.
\end{align*}
 Then for sufficiently large $m\in\mbN$ and $x\in(\e^{-2(m+1)},\e^{-2m}]$ 
\begin{align*}
\vf_1(0) - \vf_1(x) &\ge \frac{m \e^{-3m}}{3\e^4} \ge \frac{\left(\e^{-2m}\right)^{3/2}}{3\e^4} \ge \frac{x^{3/2}}{3\e^4}.
\end{align*}
On the other hand, for any $k\in\mbN$ and $x_m=\e^{-km}, m\in\mbN,$ 
\[
\vf_1(0)-\vf_1(x_m) \ge  \frac{(2k-1)m\e^{(2k-1)m} x_m^2}{3}   = \frac{2k-1}{3k} x_m^{\frac{1}{k}} \log x_m^{-1}.
\]
 
\end{example}

\begin{theorem}
\label{th:flow.existence}
Suppose $\vf\in \Phi_*$ and $a$ is measurable and of linear growth. Then the Harris flow $X^{\vf, a}$ exists and is unique in distribution.
\end{theorem}

\begin{theorem}
\label{th:flow.coalescence}
Suppose that $\vf\in \Phi_\alpha$ and $a \in A_\beta$ with $\beta - \alpha > -  1, \alpha< 2.$ Then for any $C, t \in \mbR_+$ 
\begin{equation}
\label{eq:coalescing.flow}
\E \sharp \left\{ X^{\vf, a}_{0,t}(u) \mid u \in [-C,C] \right\} < \infty.
\end{equation} 
\end{theorem}

Either result is essentially an extension of the corresponding theorem in \cite{Ha84Coalescing}. Still,  original proofs need to be modified to accommodate the presence of non-zero drift, so we present a brief description of the necessary changes in Appendix A.

\section{Harris flows as solutions to SDEs} \label{sec:sdes}

This section describes the approach of \cite{WarWa04Spectra} that provides the representation of $X^{\vf,a}$ as a solution to a SDE w.r.t a cylindrical Wiener process. 

Let $H_\vf$ be the separable Hilbert space obtained as the completion of
\[
\mathrm{span}\Big\{\sum_{k=\ov{1,n}} a_k \vf(x_k -\cdot), a_k, x_k \in \mbR, k=\ov{1,n}, n\in\mbN \Big\}
\]
w.r.t the inner product $(\vf(x-\cdot), \vf(y-\cdot))_{H_\vf} = \vf(x-y).$ 
Define
\[
M_{s,t}(x) = X^{\vf, a}_{s,t}(x) - \int_s^t a(X^{\vf,a}_{s,r}(x)) dr, \quad 0\le s\le t, x\in\mbR. 
\]
\begin{proposition}[{\cite[p. 351 + p. 356]{WarWa04Spectra}}]
\label{prop:2.1.1}
Assume that $\vf \in C(\mbR)$ and $a$ is measurable and of linear growth. Then there exists a standard cylindrical Wiener process $W$ on $H_\vf$ such that
\[
\left(W_t - W_s, \vf(x-\cdot)\right)_{H_\vf} = L_2-\lim_{n\to\infty} \sum_{k=\ov{0,n}} \left(M_{s + \frac{k}{n}(t-s),s + \frac{k+1}{n}(t-s)}(x) -x \right). 
\]
\end{proposition}

We denote such $W$ as $\cW(X^{\vf,a}).$

Assume that $e_n, n\in\mbN,$ is an orthonormal basis in $H_\vf.$ Then given $w^n_{\cdot} = (W_{\cdot}, e_n)_{H_\vf}, n\in\mbN,$
\begin{equation*}
\label{eq:2.1.add.1}
 \int_s^t \sigma(X^{\vf, a}_{s,r}(x)) dW_r = \sum_{n\ge 1} \int_s^t e_n(X^{\vf, a}_{s,r}(x)) dw^n_r = \int_s^t W\left(X^{\vf,a}_{s,r}(x), dr\right), 
\end{equation*}
where the last integral is understood in the sense of~\cite{Ku90Stochastic} and $\sigma(x), x\in\mbR,$ are  Hilbert-Schmidt operators from $H_\vf$ to $\mbR:$
\begin{equation*}
\sigma(x)(h) = \sum_{n\ge 1} (e_n, h)_{H_\vf} e_n(x), \quad h \in H_\vf.
\end{equation*}
\begin{proposition}[{\cite[p. 356]{WarWa04Spectra}}]
\label{prop:2.1.2}
Under the same assumptions as in Proposition \ref{prop:2.1.1} 
for all $x\in\mbR$ and $s\ge 0$ with probability $1$
\begin{equation}
\label{eq:2.1.add.main}
X^{\vf, a}_{s,t}(x) = x + \int_s^t a(X^{\vf, a}_{s,r}(x)) dr +  \int_s^t \sigma(X^{\vf, a}_{s,r}(x)) dW_r, \quad t\ge s.
\end{equation}
\end{proposition}
Here
\begin{align*}
\|\sigma(x)\|_{HS} &= 1, \\
\|\sigma(x) - \sigma(y)\|^2_{HS} &= 2(1 -\vf(x-y)), \quad x,y\in\mbR, 
\end{align*}
where $\|\cdot\|_{HS}$ is the corresponding Hilbert-Schmidt norm.

Both propositions are formulated in~\cite{WarWa04Spectra} without proofs and in the case of $a=0,$ so we need to justify them for nontrivial drift. Sketches of the corresponding proofs are presented in Appendix B.

 We denote the space of Holder continuous functions of order $\beta$ on $\mbR$ by $H_\beta(\mbR)$ henceforth. 
\begin{theorem}[{\cite[p. 356]{WarWa04Spectra}}]
\label{th:2.2.th.3}
 If $\sqrt{1-\vf}\in H_\beta(\mbR), \beta \ge \frac{1}{2}, a\in Lip(\mbR),$ then for every $x\in\mbR$ and $s\ge 0$ the process $X^{\vf, a}_{s,\cdot}(x)$ is the unique strong solution of \eqref{eq:2.1.add.main}.
\end{theorem}



\section{The splitting scheme and the example of the Brownian web}
 Only $a\in Lip(\mbR)$ is considered hereinafter unless stated otherwise explicitly. $T >0$ is fixed. 
We define the composition of $m\in\mbN$ functions by
\[
\mathop{\circ}_{j=\ov{1,m}} f_k = \Big( \mathop{\circ}_{j=\ov{2,m}} f_k \Big) \circ f_1.
\]  

Consider a sequence $\mathcal{T}=(\{t^n_{j} \mid j=\ov{0,N^n}\})_{ n\in\mbN}$ of partitions of $[0,T]:$ 
\begin{align*}
0 &= t^n_{0} < \ldots  t^n_{N^n} = T, \quad n\in\mbN,
\end{align*}
and set
\begin{align*}
I^n_{k} &= [t^n_{k},t^n_{k+1}), \quad k = \ov{0,N^n-2}, \\
I^n_{N^n-1} &= [t^n_{N^n-1}, T], \\
\delta^n_{k} &= t^n_{k+1} -t^n_{k}, \quad k=\ov{0,N^n-1}, \\
\delta_n &= \max_{k=\ov{0,N^n-1}} \delta^n_{k},
\end{align*}
with $\delta_n \to 0, n\to\infty.$

 Let $F_t(x), t\ge 0,$ be the solution to 
\begin{align*}
\frac{dF_t(x)}{dt} &= a(F_t(x)), \notag \\
F_0(x) &= x,  
\end{align*}
for $x\in\mbR.$

For all $ n\in\mbN$ set $y^n_{0-}(x) = x$ and define processes $(u^n, y^n)\in D([0,T],\mbR^2)$ such that  for $k =\ov{0,N^n-1}$  and $t\in  I^n_k$ (cf. \cite[Eq. 2.1]{DoVov18Arratia} for the Brownian web)
\begin{align}
\label{eq:splitting.flows}
u^n_t(x) &= F_{t - t^n_{k}}\left(y^n_{t^n_{k}-}(x)\right)  =   F_{t -t^n_{k}} \circ \Big[\mathop{\circ}_{j=\ov{0,k-1}} \Big( X^{\vf,0}_{t^n_{j}, t^n_{j+1}} \circ F_{t^n_{j+1}-t^n_{j}} \Big) \Big](x), \notag \\
y^n_t(x) &= X^{\vf,0}_{t^n_{k},t}\left(u^n_{t^{n}_{k+1} -}(x)\right) \notag \\
&\qquad = X^{\vf,0}_{t_{n_k}, t} \circ F_{t^{n}_{k+1} -t^{n}_{k}} \circ \Big[ \mathop{\circ}_{j=\ov{0,k-1}} \Big( X^{\vf,0}_{t^{n}_{j}, t^{n}_{j+1}} \circ F_{t^{n}_{j+1}-t^{n}_{j}} \Big) \Big](x).
\end{align} 

If $\vf\in \Phi_*$ and $\sqrt{1-\vf}\in H_\beta(\mbR), \beta \ge \frac{1}{2},$ \eqref{eq:splitting.flows} is equivalent in distribution to the following system. Let $W=\cW(X^{\vf,a}).$ 
Formally define 
\begin{align}
\label{eq:splitting.sdes}
u^n_t(x) &= y^n_{t^{n}_{k}-}(x) + \int_{t^{n}_{k}}^t a(u^n_s(x)) ds, \notag  \\
y^n_t(x) &= u^n_{t^{n}_{k+1} -}(x) + \int_{t^{n}_{k}}^t \sigma(y^n_s(x)) dW_s, \notag \\
y^n_{0-}(x) &= x, \notag \\
& \quad t \in I^n_{k}, k =\ov{0,N^n-1}, n\in\mbN.
\end{align} 

Borrowing notation from \cite{BenGloRas92Approximations} and setting
\begin{align*}
d^n_t &= \max\left\{t^{n}_{k} \mid t^{n}_{k} \le t\right\}, \\
\ov{d}^n_t &= \min\left\{t^{n}_{k} \mid t^{n}_{k} > t\right\}, \quad t \in [0,T), \\
d^n_T &= t^{n}_{N^n-1}, \ \ov d^n_T =T,
\end{align*}
we can rewrite  both~\eqref{eq:splitting.flows} and \eqref{eq:splitting.sdes} as 
\begin{align*}
u^n_t(x) &= x +  \int_{0}^t a(u^n_s(x)) ds + w^n_{d^n_t}(x), \notag \\
y^n_t(x) &= x +  \int_{0}^{\ov{d}^n_t} a(u^n_s(x)) ds + w^n_{t}(x), \notag \\
&\quad t\in [0,T], n\in \mbN,
\end{align*}
where  $w^n(x)$ are standard Wiener processes. In the case of~\eqref{eq:splitting.sdes} 
\[
w^n_t(x) = \int_0^{t} \sigma(y^n_s(x)) dW_s, \quad t\in [0,T], n\in \mbN.
\]

The collection $(u^n_t(\cdot),y^n_t(\cdot)), t\in[0,T],$ can be considered as a $D^\uparrow(\mbR^2)-$valued random process in either case.
\begin{definition}
 We denote  $((u^n,y^n))_{n\in\mbN}$ by $\mathrm{Spl}(X^{\vf,0}; a; \mathcal{T})$  in the case of~\eqref{eq:splitting.flows} and by $\wt{\mathrm{Spl}}(X^{\vf,a}; a; \mathcal{T})$ in the case of~\eqref{eq:splitting.sdes}, respectively. 
\end{definition}

We summarize the discussion of the well-posedness of~\eqref{eq:splitting.flows} and \eqref{eq:splitting.sdes} as follows.

\begin{proposition}
\label{prop:splitting}
\begin{enumerate}
\item If $\vf\in \Phi_*, a\in Lip(\mbR),$ then $\mathrm{Spl}(X^{\vf,0}; a; \mathcal{T})$ is unique in distribution. 
\item  If $\sqrt{1-\vf}\in H_\beta(\mbR), \beta \ge \frac{1}{2},$ additionally, then for any $x\in\mbR$ and $n\in\mbN$ the pair $(y^n_\cdot(x), u^n_\cdot(x))$ from $\wt{\mathrm{Spl}}(X^{\vf,a}; a; \mathcal{T})$ is the  unique strong $(\cF^{X^{\vf,a}}_{0,s})_{s\in[0,T]}-$\-adapted solution to the system  \eqref{eq:splitting.sdes}; moreover,  $\mathrm{Spl}(X^{\vf,0}; a; \mathcal{T})$ and \\ $\wt{\mathrm{Spl}}(X^{\vf,a}; a; \mathcal{T})$ are identically distributed.
\end{enumerate}
\end{proposition} 

The proof of the following lemma uses reasoning similar to that in \cite[pp. 171--172]{KaShre91Brownian} and is therefore omitted.
\begin{lemma}
\label{lem:hol.2}
Assume that  $((u^n,y^n))_{n\in\mbN}=\wt{\mathrm{Spl}}(X^{\vf,a}; a; \mathcal{T}),$  $W=\cW(X^{\vf,a})$ and  $\vf\in \Phi_*,$  $\sqrt{1-\vf}\in H_\beta(\mbR)$  for some $\beta\ge \frac{1}{2}.$ For any  $x\in\mbR$ and $n\in\mbN$ there exists, possibly on an extension of the original probability space, a standard Wiener process $b^n(x)$ such that
\begin{align*}
\int_0^{t} \left(\sigma\left(X^{\vf,a}_{0,s}(x)\right) - \sigma\left(y^n_s(x)\right)\right) dW_s &= 2^{1/2}\int_0^t \left( 1 -\vf\left(X^{\vf,a}_{0,s}(x) - y^n_s(x)\right)\right)^{1/2} db^n_s(x), \\ 
& t\in [0,T].
\end{align*}
\end{lemma}


Define 
\begin{align*}
r^n_t(x) &= - \int_t^{\ov{d}^n_t} a(u^n_s(x)) ds, \notag \\
l^n_t(x) &=  w^n_t(x) - w^n_{d^n_t}(x), \notag  \\
t &\in [0,T], x \in\mbR, n\in\mbN, 
\end{align*} 
so 
\begin{align}
\label{eq:notation.1}
y^n(x) - u^n(x) = l^n(x) + r^n(x). 
\end{align}

\begin{lemma}
\label{lem:hol.3}
Let  $((u^n,y^n))_{n\in\mbN}=\mathrm{Spl}(X^{\vf,0}; a; \mathcal{T}).$
\begin{enumerate}
\item For any $x\in\mbR$ and $p\ge 1$ there exists  $C = C(p, x, T)$ such that
\begin{align*}
\sup_{t\in[0,T]}|u^n_t(x)|^p &\le C \Big(1 + \sup_{t\in[0,T]}\left| w^n_t(x) \right|\Big)^p, \quad n\in\mbN, \\
\sup_{n\in\mbN} \E \sup_{t\in[0,T]}  |u^n_t(x)|^p &\le C.
\end{align*}
\item For $p\ge 2$  there exists $C_1 = C_1(p)$ such that
\begin{align*}
 \E \sup_{t\in[0,T]} |r^n_t(x)|^p &\le T C_1 \E \Big( \int_{0}^T \left( 1+  \left| u^n_s(x) \right| \right)^2 ds \Big)^{p/2}  \, (\delta_n)^{p/2-1} , \quad x\in\mbR, n\in\mbN, \\
 \sup_{x\in\mbR} \E \sup_{t\in[0,T]} |l^n_t(x)|^p &\le T C_1  \, (\delta_n)^{p/2-1}, \quad n\in\mbN.  
\end{align*}
\item For any $x\in\mbR$ and $p\ge 2$ there exists  $C_2 = C_2(x,T)$ such that
\begin{align*}
\sup_{n\in\mbN} \E \sup_{t\in[0,T]} |r^n_t(x)|^p &\le C_2, \\
\sup_{n\in\mbN} \E \sup_{t\in[0,T]} |l^n_t(x)|^p &\le C_2. 
\end{align*}
\item There exists $C_3=C_3(T)$
\begin{align*}
\sup_{x\in\mbR} \E \sup_{t\in[0,T]} \left( l^n_t(x)\right)^2 \le C_3 \delta_n \log \delta_n^{-1}, \quad n\in\mbN. 
\end{align*}
\item For any $x\in\mbR$ there exists $C_4=C_4(x,T)$
\begin{align*}
\E \sup_{t\in[0,T]} \left( r^n_t(x)\right)^2 \le C_4 \delta_n, \quad n\in\mbN. 
\end{align*}
\end{enumerate}
\end{lemma}
\proof We drop the argument $x$ which is assumed to be fixed.

(1) The inequality follows trivially by the Gronwall lemma.  

(2) We have for some $\wt C$ 
\begin{align*}
\E \sup_{t\in[0,T]} \left| r^n_t \right|^p &\le \sum_{k=\ov{0, N^n-1}} \E \sup_{t\in I^n_{k}} \Big| \int_t^{\ov{d}^n_t} a\left( u^n_s\right) ds \Big|^p \\
&\le    \wt C \E \Big( \int_{0}^T \left( 1+  \left| u^n_s \right| \right)^2 ds \Big)^{p/2} \sum_{k=\ov{0, N^n-1}} (\delta^n_{k})^{p/2}, \\
\E \sup_{t\in[0,T]} \left| l^n_t \right|^p &\le \sum_{k=\ov{0, N^n-1}} \E \sup_{t\in I^n_{k}} \left| w^n_{t} - w^n_{d^n_t} \right|^p \\
&\le  \wt C \sum_{k=\ov{0, N^n-1}} (\delta^n_k)^{p/2}.
\end{align*}

(3) is a corollary of (1) and (2).

(4) Set 
\begin{align*}
\xi^n_{k} &= \sup_{s \in I^n_{k} } \left(w^n_{s} - w^n_{t^{n}_{k}} \right)^2, \quad k =\ov{0,N^n-1}.
\end{align*}
Then we need to estimate
\begin{align*}
\E \sup_{t \in [0,T]} \left(w^n_{t} - w^n_{d^n_t} \right)^2 &= \int_0^\infty \Big[1 - \Prob\Big( \max_{k=\ov{0,N^n-1}} \xi^n_{k}  \le u \Big)\Big] du.
\end{align*}
Since
\[
\Prob\Big( \max_{k=\ov{0,N^n-1}} \xi^n_{k}  \le u \Big) = \E  \Prob\Big(   \xi^n_{N^n-1} \le u | \cF^{X^{\vf,0}}_{t^{n}_{N^n-1}} \Big) \prod_{k=\ov{0,N^n-2}} \1\left[ \xi^n_{k} \le u\right]
\]
and, for a standard Wiener process $w$ and $k=\ov{0,N^n-1},$  
\begin{align*}
\Prob\left(  \xi^n_{k} \le u | \cF^{X^{\vf,0}}_{t^{n}_{k}} \right) &\ge 1 - 2\Prob\Big( \sup_{s \in I^n_{k}} w_s \ge u^{1/2} \Big) \\ 
&= 1 - \frac{4}{(2\pi)^{1/2}} \int_{\sqrt{\frac{u}{\delta^n_{k}}}}^\infty \e^{-\frac{v^2}{2}} dv,
\end{align*}
we get, iterating conditioning, using the inequality
\[
\prod_{k=\ov{1,m}} (1 -x_k) \ge 1 - \sum_{k=\ov{1,m}} x_k, \quad x_k \in [0,1], k =\ov{1,m}, m\in\mbN,
\]
and standard estimates for the Gaussian distribution, that given $\alpha=2 \delta_n \log \frac{1}{\delta_n}$ we have for sufficiently large $n$ and some absolute constants $K_1, K_2$ 
\begin{align*}
\E \sup_{t \in [0,T]} \left(w^n_{t} - w^n_{d^n_t} \right)^2 &\le \alpha  + K_1 \int_{\alpha}^\infty \Big( \sum_{k=\ov{0,N^n-1}}  \int_{\sqrt{\frac{u}{\delta^n_{n}}}}^\infty \e^{-\frac{v^2}{2}}dv \Big) du \\
&\le \alpha  + K_1 \sum_{k=\ov{0,N^n-1}}  \int_{\alpha}^{\infty }  \e^{-\frac{u}{2\delta^n_k}} \sqrt{\frac{\delta^n_{k}}{u}} du \\
&\le \alpha  + \frac{K_2}{ \sqrt{\alpha} } \sum_{k=\ov{0,N^n-1}} (\delta^n_{k})^{3/2} \cdot \int_{ \frac{\alpha}{2\delta^{n}} }^{\infty } \e^{-u} du   \\
&= \alpha +  K_2 T  \e^{-\frac{\alpha}{2\delta_n}} \sqrt{\frac{\delta_n}{\alpha}}  \\
&= 2 \delta_n \Big( \log \frac{1}{\delta_n} + \frac{K_2 T}{\sqrt{2\log \frac{1}{\delta_n}}}\Big). 
\end{align*} 

(5) For some $\wt C$
\begin{align*}
\E \sup_{t\in[0,T]} \left(r^n_t\right)^2 &\le \wt C \E \sup_{k=\ov{0,N^n-1}} \delta^n_{k} \int_{I^n_{k}} \left(1 + \left( u^n_s \right)^2  \right) ds \\
&\le \wt C \delta_{n} \E \int_0^T \left(1 + \left( u^n_s \right)^2  \right) ds,
\end{align*} so the application of (1) yields the desired estimate. \qed

Let $\cM_p(\mbR)$ be the metric space of probability measures on $\mbR$ with finite $p-$th moment and let $W_p$ be the corresponding Wasserstein distance~\cite{VillaniTopics03}. 
For fixed $p,$ we define the Wasserstein distance between probability measures $L_1, L_2$ on $\cM_p(\mbR)$ 
 as
$$
W_{1, p}(L_1, L_2)=\inf\E \ \! W_p(\mu^\prime, \mu^{\prime\prime}),
$$
where the infinum is taken over the set of pairs of $\cM_p(\mbR)-$valued random elements $\mu_1, \mu_2$ satisfying $\Law(\mu_k) = L_k, k=1,2.$

Define random pushforward measures 
\begin{align}
\label{add:10.1}
\mu_t &= \Leb_{[0, 1]}\circ\left(X^{\vf,a}_t\right)^{-1}, \notag \\
\mu^{n}_t &= \Leb_{[0, 1]}\circ\left(y^n_t\right)^{-1}, \quad t\in[0,T], n\in\mbN,
\end{align}
and their distributions as measures on $\cM_p(\mbR)$
\begin{align*}
L_t &= \Law(\mu_t), \\
L^{n}_t &= \Law\left(\mu^{n}_t\right), \quad t\in[0,T], n\in\mbN.
\end{align*}

The rest of the section describes the example of splitting for the Brownian web. Let $B$ be a Brownian web. The corresponding counterpart with drift $B^a$ is defined and constructed in \cite[Chapter 7]{Do07MeasureEng} as a family $\{B^a_\cdot(x)\mid x \in\mbR\}$ of coalescing semimartingales.  One defines the associated splitting $\mathrm{Spl}(B; a; \mathcal{T})$ via \eqref{eq:splitting.flows} by replacing $X^{\vf,0}$ with $B.$
It is worth noting that  the limit in Proposition \ref{prop:2.1.1} does not exist due to \cite[Proposition 1.5]{DoVov18Arratia}. 
\begin{theorem}[{\cite[Theorem 4.1]{DoVov18Arratia}}]
\label{th:split.1}
Assume $a\in L_\infty(\mbR).$ For any $m\in\mbN$ and any $x_1,\ldots,x_m\in \mbR$
\[
\left(y^n_\cdot(x_1),\ldots, y^n_\cdot(x_m)\right) \Rightarrow \left( B^a_\cdot(x_1), \ldots, B^a_\cdot(x_m)\right), \quad n\to\infty, 
\]
in $D([0,T],\mbR^m).$
\end{theorem}

Set $L_t = \Law(\Leb_{[0,1]}\circ (B^a_t)^{-1}).$ 
\begin{theorem}[{\cite[Theorem 2.1]{DorVov20ApproximationsEng}}]
\label{th:split.2}
Assume that the sequence $\left\{n\delta_n \right\}_{n\in\mbN}$ is bounded by $K$ and $a\in L_\infty(\mbR).$ Then for every $p\geq 2$ there
exist $C= C(p, K, T)>0$  such that 
\[
W_{1,p}(L_t, L^{n}_t)\le \frac{C}{\left(\log\log \delta_n^{-1}\right)^{1/p}}, \quad n\in\mbN.
\]
\end{theorem}


\begin{remark}
The formulation of \cite[Theorem 2.1]{DorVov20ApproximationsEng} is erroneously missing the second logarithm due to a calculational error in the end of the proof.
\end{remark}

\section{Weak convergence}
Let $((y^n,u^n))_{n\in\mbN}=\mathrm{Spl}(X^{\vf,0}; a; \mathcal{T})$ for some $\mathcal{T}.$ The main result of this section is the following theorem.
\begin{theorem}
\label{th:con.1}
Assume that $\vf \in\Phi_*$ and $a\in Lip(\mbR).$ For any $m\in\mbN$ and any $x_1,\ldots,x_m\in \mbR$
\[
\left(y^n_\cdot(x_1),\ldots, y^n_\cdot(x_m)\right) \Rightarrow \left( X^{\vf,a}_{0,\cdot}(x_1), \ldots, X^{\vf,a}_{0,\cdot}(x_m)\right), \quad n\to\infty, 
\]
in $D([0,T],\mbR^m).$
\end{theorem}
The proof of Theorem \ref{th:con.1} is split into a series of lemmas. 

Recall  
\[
w^n_t(x) = y^n_t(x) - x - \int_0^{\ov{d}^n_t} a\left( u^n_s(x)\right) ds, \quad t \in[0,T], x\in\mbR,  n\in\mbN,
\]
to be standard Wiener processes. 

We denote the modulus of continuity by $\omega$ and the Lipschitz constant for $a$ by $C_a.$

\begin{remark}
Proceeding exactly as in the proof of~\cite[Proposition 2.2]{DoVov18Arratia}, one can show that 
for $u_1,u_2\in\mbR$ 
\[
\langle w^n(u_1), w^n(u_2) \rangle_t = \int_0^t \vf\left( y^n_s(u_1) - y^n_s(u_2)\right) ds, \quad t\in[0,T], n\in\mbN. 
\]
\end{remark}

Denote
\begin{align*}
x &= (x_1,\ldots, x_m), \\
y^n(x) &= (y^n_\cdot(x_1),\ldots, y^n_\cdot(x_m)), \\
u^n(x) &= (u^n_\cdot(x_1),\ldots, u^n_\cdot(x_m)), \\
& \qquad n\in\mbN.
\end{align*} 
\begin{lemma}
\label{lem:con.2}
The sequence $(u^n(x), y^n(x))_{n\in\mbN}$ is weakly relatively compact in $D([0,T],\mbR^{2m}).$
\end{lemma}
\proof By \cite[Theorem 15.2]{Bi68Convergence}, it is sufficient to prove that for any $\ve>0$ and any $k=\ov{1,m}$ 
\begin{equation}
\label{eq:add.10}
\lim_{\vk\to 0+} \limsup_{n\to \infty} \Prob\left( \max\left\{\wh\omega(y^n(x_k),\vk), \wh\omega(u^n(x_k),\vk)\right\} \ge \ve \right)=0,
\end{equation}
where for $f\in D([0,T],\mbR)$
\[
\wh\omega(f, \vk) = \inf_{\begin{subarray}{c}0=t_0 < t_1 < \ldots < t_r = T,\\ t_i - t_{i-1}>\vk, i=\ov{1,r}, \\
r \in \mbN \end{subarray} } \max_{i=\ov{1,r}} \sup_{s_1, s_2 \in [t_{i-1},t_i)} \left(f_{s_1} - f_{s_2}\right). 
\]
For fixed $k$ and $\ve$ consider $s_1, s_2$ such that $0 < s_2 -s_1 < \vk.$ We drop argument $x$ to simplify notation. Since
\begin{align*}
 y^n_{s_2} - y^n_{s_1} &= \int_{\ov{d}^n_{s_1}}^{\ov{d}^n_{s_2}} a(u^n_r) dr + w^n_{s_2} - w^n_{s_1}, \\
 u^n_{s_2} - u^n_{s_1} &= \int_{s_1}^{s_2} a(u^n_r) dr + w^n_{d^n_{s_2}} - w^n_{d^n_{s_1}},
\end{align*}
we have for some $C>0$
\begin{align*}
\omega\left(y^n,\vk\right) &\le \omega\left(w^n, \vk\right) + C\Big(1 + \sup_{t\in[0,T]} \left|u^n_t\right|\Big) (\vk + 2\delta_n), \\
\omega\left(u^n,\vk\right) &\le \omega\left(w^n, \vk + 2\delta_n\right) + C\Big(1 + \sup_{t\in[0,T]} \left|u^n_t\right|\Big) \vk.
\end{align*}
Let $w$ be a standard Wiener process. By Lemma~\ref{lem:hol.3}, for some fixed constant $K$ 
\begin{multline}
\label{eq:add.13}
\Prob\left( \max\left\{\omega(y^n(x_k),\vk), \omega(u^n(x_k),\vk)\right\} \ge \ve \right) \le \\
\Prob\left( \omega(w, \vk + 2\delta_n) \ge \frac{\ve}{2} \right)  + \Prob\Big( \Big(1+\sup_{t\in[0,T]} |w_t|\Big) (\vk +\delta_n)  \ge \frac{\ve}{K}\Big)
\end{multline}

Since $\wh\omega(\cdot, \vk) \le \omega(\cdot, 2\vk)$ for $\vk \ll 1,$ \eqref{eq:add.10} follows. \qed


\begin{lemma}
\label{lem:con.4} 
For any weak limit $\xi=(\xi_1,\ldots, \xi_m)$ of the sequence $(y^n(x))_{n\in\mbN}$  and  for any pair $i, j\in\{1,\ldots,m\}, i\not=j,$
\begin{align*}
\Prob\Big(\exists t\in [0; T] \ \xi_{it}=\xi_{jt}\mbox{\ and\ } \sup_{s\in[t; T]}|\xi_{is} - \xi_{js}|>0 \Big) = 0.
\end{align*} 
\end{lemma}
\proof We adopt the idea from~\cite[Proposition 3.8]{DoVov18Arratia}, referring to the aforementioned proof for those calculations that are shared between the proofs. Define
\begin{align*}
D^+([0,T],\mbR) &= \Big\{f\in D([0,T],\mbR)\mid\inf_{r\in[0,T]}f_r\ge 0 \Big\}, \\
\Gamma^\vk_\ve &= \Big\{f\in D^+([0,T],\mbR)\mid\exists t\in[0,T] \colon \ f_t<\ve, \int^T_t f_r dr>\vk\Big\},
\\
\Gamma^\vk &= \Big\{f\in D^+([0,T],\mbR)\mid\exists t\in[0,T]  \colon f_t=0, \int^T_t f_r dr>\vk\Big\}, \\
& \vk, \ve > 0.
\end{align*}
Assume that $i, j$ are fixed and that $x_j > x_i.$ It is sufficient to show that for any $\vk$ 
\begin{equation}
\label{eq:add.11}
\liminf_{\ve\to 0+} \liminf_{n\to\infty} \Prob\left( \Delta y^n \in \Gamma^\vk_\ve\right) = 0.
\end{equation}
where $\Delta y^n=y^n(x_j) - y^n(x_i).$  Put 
\[
\omega^n = \omega(w^n(x_i), \delta_n) + \omega(w^n(x_j), \delta_n).
\]
One can show that for $\ve \ll \vk$
\begin{align*}
\Prob\left( \Delta y^n \in \Gamma^\vk_\ve\right)  &= \Prob\Big(\exists t\colon \Delta y^n_{t}<\ve, \int^T_t \Delta y^n_{r}dr>\vk\Big)  \\
&\leq 
\Prob\left(\omega^{n}\geq \ve \right) +  \sum_{l= \ov{0,N^n-1}}
\Prob\Big( \inf_{r\in[0, t^{n}_{l-1})} \Delta y^n_r\geq\ve; \inf_{r\in[t^{n}_{l-1}, t^{n}_{l})} \Delta y^n_r\le \ve;  \\ 
& \qquad\qquad \Delta y^n_{t^{n}_{l}-}\leq2\ve; \int_{t^{n}_{l}}^T \Delta y^n_r dr \geq \frac{\vk}{2} \Big) \\
&= \Prob\left( \omega^{n}\geq \ve \right) + I^n,
\end{align*} 
where $t^{n}_{-1}=0.$ Since 
\[
\lim_{n\to\infty}\Prob\left\{\omega^n \geq \ve\right\} =0,
\]
we consider only the sum $I^n.$ For $t \in [t^{n}_{p},t^{n}_{p+1}), p \ge  l,$ %
\begin{align}
\label{eq:add.16}
\E \Big( \Delta y^n_t \vert \cF^{X^{\vf,0}}_{t^{n}_{l}}\Big) &= \E \Big( \E \Big(  \Delta y^n_t \vert \cF^{X^{\vf,0}}_{t^{n}_{p}}\Big) \vert \cF^{X^{\vf,0}}_{t^{n}_{l}}\Big) \notag \\
&= \E \Big( \Delta y^n_{t^{n}_{p}}  \vert \cF^{X^{\vf,0}}_{t^{n}_{l}}\Big) \notag  \\
&\le \e^{C_a \delta^n_{p}} \E \Big( \Delta y^n_{t^{n}_{p}-}  \vert \cF^{X^{\vf,0}}_{t^{n}_{l}}\Big) \notag  \\
&\le \ldots \notag  \\
&\le \e^{C_a(T-t^{n}_{l})} \Delta y^n_{t^{n}_{l}-},
\end{align}
so
\begin{align}
\label{eq:add.16.1}
\Prob\Big( \int_{t^{n}_{l}}^T \Delta y^{n}_r dr \geq \frac{\vk}{2} \Big\vert \cF^{X^{\vf,0}}_{t^{n}_{l}} \Big) &\leq 
\frac{2}{\vk} \E\Big(\int_{t^{n}_{l}}^T \Delta y^{n}_r dr \Big\vert \cF^{X^{\vf,0}}_{t^{n}_{l}} \Big) \notag \\
&\leq \frac{2 \e^{C_a(T-t^{n}_{l})} (T-t^{n}_{l})}{\vk} \Delta y^n_{t^{n}_{l}-}.
\end{align}
Thus
\begin{align*}
I^n &\le \frac{4\e^{C_aT}T}{\vk} \ve,
\end{align*}
which yields \eqref{eq:add.11}. \qed

\begin{lemma}
\label{lem:con.3}
Let $\mcC_m$ be the set of elements of $C(\mbR_+, \mbR^m)$ whose coordinates merge after a meeting. 
Any weak limit $\xi$ of the sequence $(y^n(x))_{n\in\mbN}$ is a $\mcC_m-$solution in the sense of Definition \ref{defn:C.solution} in Appendix A to the martingale problem on $\mbR^m$ for the operator
\[
\cA_m = \frac{1}{2}\sum_{k,j=\ov{1,m}} \vf(x_k-x_j) \frac{\pt^2}{\pt x_k \pt x_j} + \sum_{k=\ov{1,m}} a(x_k) \frac{\pt }{\pt x_k}.
\] 
\end{lemma}
\proof W.l.o.g. we can suppose that $y^n(x) \Rightarrow \xi, n\to\infty.$  Using \cite[Theorem 15.5]{Bi68Convergence} and \eqref{eq:add.13} we can check that $\xi$ is continuous a.s.. Recalling~\eqref{eq:notation.1},  define
\begin{align*}
v^n(x_k) &= y^n(x_k) - r^n(x_k), \quad k=\ov{1,m}, \\
v^n(x) &= (v^n(x_1),\ldots, v^n(x_m)),  \\ 
u^n(x) &= (u^n(x_1),\ldots, u^n(x_m)), \quad n\in\mbN. 
\end{align*}
Lemma \ref{lem:hol.3} implies
\begin{align*}
 \left(r^n(x_1),\ldots, r^n(x_m)\right) &\to 0,\quad  n\to\infty, \notag \\
 \left(l^n(x_1),\ldots, l^n(x_m)\right) &\to 0, \quad n\to\infty, 
\end{align*}
in probability in $D([0,T],\mbR^m)$  in the uniform metric and therefore in the $J_1$ topology. Thus, 
\begin{align}
\label{eq:add.15}
 y^n(x) - u^n(x) &\to 0,\quad  n\to\infty, \notag \\
 y^n(x) - v^n(x) &\to 0, \quad n\to\infty, 
\end{align} 
in probability and therefore 
\begin{align*}
v^n(x)  &\Rightarrow \xi, \quad n\to\infty, \\
u^n(x) &\Rightarrow \xi, \quad n\to\infty,
\end{align*} 
in $D([0,T],\mbR^m).$ Proceeding as in Lemma~\ref{lem:con.2}, we can check that the sequence 
 \[
 \left((y^n(x),u^n(x), v^n(x))\right)_{n\in\mbN}
 \] is weakly relatively compact in $D([0,T],\mbR^{3m}).$ Hence, by the Skorokhod representation theorem  and \eqref{eq:add.15} we can assume w.l.o.g that 
\begin{align}
\label{eq:add.160}
(y^n(x),u^n(x),v^n(x)) &\Rightarrow (\xi, \xi, \xi), \quad n\to\infty, 
\end{align}
in  $D([0,T],\mbR^{3m}).$
 
By the It\^o   lemma and Proposition~\ref{prop:2.1.2}, for any bounded $f\in C^2(\mbR^m)$ 
\begin{align*}
f(v^n_t(x)) &= f(x) + \int_0^t \sum_{k=\ov{1,m}} a(u^n_s(x_k)) \frac{\pt }{\pt x_k} f(v^n_s(x)) ds \\
 & \qquad + \int_0^t \sum_{j,k=\ov{1,m}} \vf\left(y^n_s(x_k) - y^n_s(x_j)\right) \frac{\pt^2 }{\pt x_k \pt x_j} f(v^n_s(x)) ds \\
 & \qquad + \int_0^t \sum_{k=\ov{1,m}} \sigma(y^n_s(x_k)) \frac{\pt }{\pt x_k} f(v^n_s(x)) dW_s,
\end{align*}
where $W=\cW(X^{\vf,0}).$

Assume $g \in C(\mbR^{Mm})\cap L_\infty(\mbR^{Mm})$ for some $M\in\mbN.$ Then for arbitrary $s,t, s \le t$ and $s_1, \ldots, s_M \le s$
\begin{align*}
\E g(y^n_{s_1}(x), \ldots, y^n_{s_M}(x)) \int_s^t \sum_{k=\ov{1,m}} \sigma(y^n_s(x_k)) \frac{\pt }{\pt x_k} f(v^n_s(x)) dW_s &= 0, 
\end{align*}
so the process
\begin{align*}
m^n_t(x) &= f(v^n_t(x)) - \int_0^t \sum_{k=\ov{1,m}} a(u^n_s(x_k)) \frac{\pt }{\pt x_k} f(v^n_s(x)) ds \\
& \qquad - \int_0^t \sum_{j,k=\ov{1,m}} \vf\left(y^n_s(x_k) - y^n_s(x_j)\right) \frac{\pt^2 }{\pt x_k \pt x_j} f(v^n_s(x)) ds, \\
& t\in[0,T],
\end{align*}
is a martingale w.r.t. the filtration generated by $y^n(x)$. 

Applying the Skorokhod representation theorem to~\eqref{eq:add.160}, we can assume that 
\begin{align*}
(y^n(x),u^n(x),v^n(x)) &\to (\xi, \xi, \xi), \quad n\to\infty, 
\end{align*}
a.s. in $D([0,T],\mbR^{3m}).$ Since the limit $\xi$ is continuous, the convergence is uniform. In particular,  for any $j,k =\ov{1,m}$
\[
y^n(x_k) - y^n(x_j)\to \xi_k - \xi_j, \quad n\to\infty,
\]
uniformly. Thus one can check that
\begin{align*}
& \E g(y^n_{s_1}(x), \ldots, y^n_{s_M}(x)) \left(m^n_t(x) - m^n_s(x) \right) \to  \\
& \qquad \E g(\xi_{s_1}(x), \ldots, \xi_{s_M}(x)) \Big( f(\xi_t) - f(\xi_s) - \int_s^t  \cA_m f(\xi_r) dr \Big), \quad  n\to\infty,
\end{align*}
so the process $t\mapsto f(\xi_t) - \int_0^t  \cA_m f(\xi_r) dr$ is a martingale.

By Lemma \ref{lem:con.4} 
\[
\Prob \left( \xi \in \mcC_m \right) = 1,
\]
 which concludes the proof. \qed

\begin{lemma}
\label{lem:con.5}
For any weak limit $\xi$ of the sequence $(y^n(x))_{n\in\mbN}$ 
\[
\Law\left( \xi \right) = \Law\left( \left(X^{\vf,a}_{0,\cdot}(x_1), \ldots, X^{\vf,a}_{0,\cdot}(x_m)\right) \right).
\]
\end{lemma}
\proof $\mcC_m-$solutions are unique by Proposition \ref{3.}. \qed

This finishes the proof of Theorem \ref{th:con.1}.

\begin{remark}
\label{rem:weak.1}
The splitting scheme and Theorem~\ref{th:con.1} can be extended to some classes of non-Lipschitz $a$ as follows. Assume that $a$ satisfies the one-sided Lipschitz condition: for some $C$
\[
a(x) - a(y) \le C (x-y), \quad x \ge y.
\] 
Then the unique flow $X^{\vf,a}$ exists and has finite moments of any order. For any $s\in [0,T)$ consider a SDE 
\begin{align*}
dF^{n}_{s,t}(x) &= a\left( F^{n}_{s,t}(x) \right) dt + \ve_n dw_t, \\
F^n_{s,s}(x) &=x,
\end{align*}
where $\ve_n\to 0+, n\to\infty,$ and $w$ is a Wiener process on $[0,T]$ independent of $X^{\vf,0}.$ Such SDEs have unique strong solutions. At each step of the splitting procedure, replace  $u^n$ on $I^n_k $ in~\eqref{eq:splitting.flows}  with 
\[
u^n_t = F^{n}_{t^n_k,t} \left( y^n_{t^n_k-}\right).
\]
Then the new splitting scheme $\mathrm{Spl}(X^{\vf,0}; a; \mathcal{T})$ is well defined. Moreover, Lemmas~\ref{lem:hol.3}~and~\ref{lem:con.2} follow immediately. For Lemma~\ref{lem:con.4}, note that general comparison theorems for SDEs~\cite{SzpruchZhang18Vintegrability, FjordMuschPi22zero_noise, Na73comparison} imply that $F^{n}_t$ is monotone mapping for any $t$ and $n$ so we can use the Gronwall lemma to establish analogs of \eqref{eq:add.16}~and~\eqref{eq:add.16.1} for conditional expectations w.r.t. the extended filtration
\[
\sigma\left\{ w_{u_1}, X^{\vf,0}_{u_1,u_2}, 0 \le  u_1 \le u_2 \le s\right\}, \quad 0 \le s \le T,
\] 
instead of $(\cF^{X^{\vf,0}}_{0,s})_{s\in[0,T]},$ which yields the conclusion of Lemma~\ref{lem:con.4}  for such $a.$ Lemma~\ref{lem:con.3} is also valid for such drift. This establishes the claim. See also~\cite{Yan02Euler} for the Euler-Maruyama scheme for SDEs with discontinuous coefficients and \cite{BuckSamTam22Splitting, BreGou19Analysis, CuiHong19strong, BreCuiHong19strong}  for splitting schemes for SPDEs with non-Lipschitz coefficients. 
\end{remark}

\section{Convergence of pushforward measures}
$((u^n,y^n))_{n\in\mbN}=\mathrm{Spl}(X^{\vf,0}; a; \mathcal{T})$  is considered in this section.

Let $\cR(\mbR)$ be the set of Radon measures on $\mbR.$ 

\begin{theorem}
\label{th:add.1}
Assume that $\vf\in \Phi_\alpha$ for some $\alpha<2$ and $\nu_0\in\cR(\mbR)$ is such that
\[
\forall \gamma > 0 \int_\mbR \e^{-\gamma u^2} \nu_0(du) < \infty.
\]
Define 
\begin{align*}
\nu_t &= \nu_0\circ\left(X^{\vf,a}_t\right)^{-1}, \\
\nu^{n}_t &= \nu_0\circ\left(y^n_t\right)^{-1}, \quad t\in[0,T], n\in\mbN.
\end{align*}
Then for any $t\in[0,T]$  $\nu^n_t, \nu_t \in \cR(\mbR)$ a.s. and
\begin{align*}
 \nu^n_t \Rightarrow \nu_t, \quad n\to \infty, 
\end{align*}
in $\cR(\mbR)$ under the vague topology. 
\end{theorem}

For the pushforward measures defined in~\eqref{add:10.1}, the following conclusion holds.

\begin{corollary}
\label{cor:add.2}
Assume that $\vf\in \Phi_\alpha$ for some $\alpha<2.$ Then for any $t\in[0,T]$  
\begin{align*}
 \mu^n_t \Rightarrow \mu_t, \quad n\to \infty, 
\end{align*}
in $\cR(\mbR)$ under the weak topology.  
\end{corollary}

We need the following lemma whose proof is postponed until Appendix A.

\begin{lemma}
\label{lem:add.1}
\begin{enumerate}
\item For all $n\in\mbN$
\[
\sup_{t\in [0,T]} \E \left| y^n_{t}(v) - y^n_{t}(w) \right| \le \e^{C_a T} |v-w|, \quad v,w \in \mbR.
\]
\item Assume that $\vf\in \Phi_\alpha$ for some $\alpha<2.$  Then for any $p\ge 2$ and $\ve\in (0;\frac{1}{1+\frac{(2-\alpha)p}{2}})$ there exists $C=C(p,\ve,T)>0$ such that
\[
 \E \sup_{t\in [0,T]} \left| X^{\vf,a}_{0,t}(v) - X^{\vf,a}_{0,t}(w) \right|^p  \le C |v-w|^\ve, \quad |v-w|\le 1.
\]
\end{enumerate}
\end{lemma}
{\it Proof of Theorem \ref{th:add.1}.}  Using Lemma~\ref{lem:add.1}, we can repeat the reasoning in \cite[proof of Theorem 1, pp. 87--90]{Vov18Convergence} as soon as we have proved the following two claims: for arbitrary $\vk >0$ and compactly supported $f$ there exists $M=M(\vk, f)$ such that 
\begin{align*}
\max\Big\{ \sup_{n\in\mbN}\E\Big|\int_{|v|\ge M} f\left(y^n_t(v)\right)d\nu_0(v)\Big|, \E\Big|\int_{|v|\ge M} f\left(X^{\vf,a}_{0,t}(v)\right)d\nu_0(v)\Big| \Big\} \le \vk,
\end{align*}
and for arbitrary $M>0$ 
\begin{align*} 
\max \Big\{ \sup_{n\in\mbN} \E \nu^n_t\left((-M;M]\right); \E \nu_t\left((-M;M]\right)  \Big\} < \infty. 
\end{align*}
For that, it is sufficient to show that for any $S>0$
\begin{align}
\label{eq:add.21}
\lim_{M\to\infty}  \int_{|x| \ge M}  \Prob\left( \left| X^{\vf,a}_{0,t}(x) \right| \le S\right)  d\nu_0(x) &= 0,  \\ 
\label{eq:add.22}
\lim_{M\to\infty} \sup_{n\in\mbN} \int_{|x| \ge M} \Prob\left( \left| y^{n}_{t}(x) \right| \le S\right) d\nu_0(x) &= 0.
\end{align}

Let $C_a$ be such that $|a(x)| \le C_a(1 + |x|)$ on $\mbR.$ The solution to the ODE 
\[
\frac{dg_t}{dt}= -C_a(1 + g_t)
\]
 is
\[
g_t =  \e^{-C_a t} -1 + \e^{-C_a t} g_0.
\]
Assume $x \gg S.$
Setting 
\begin{align*}
\zeta^n_{k}(x) &=  w^n_{t^{n}_{k+1}}(x) -w^n_{t^{n}_{k}}(x), \quad k=\ov{0,N^n-1}, n\in\mbN, x\in\mbR, 
\end{align*}
 define 
\begin{align*}
\eta^n_t(x) &= \e^{-C_a t^{n}_{k}}x  + w^n_t(x) -w^n_{t^{n}_{k+1}}(x) + \sum_{j=\ov{0,k}} \e^{-C_a (t^{n}_{k+1}-t^{n}_{j+1}) } \left(\e^{-C_a \delta^n_{j}} - 1 + \zeta^n_{j}  \right), \\
& t \in I^n_{k}, k =\ov{0,N^n-1}, n\in\mbR, x\in\mbR.
\end{align*}
Then on $\{ \inf_{t\in[0,T]} \eta^n_t(x) > 0 \}$ for $t\in I^n_{k}$ for some $k$
\begin{align*}
y^n_t(x) &\ge w^n_t(x) -w^n_{t^{n}_{k}}(x)  + \e^{-C_a \delta^n_{k}} - 1 + \e^{-C_a \delta^n_{k}} y^n_{t^{n}_{k}-} \\
&\ge \eta^n_t(x).
\end{align*}
Thus
\[
\Prob\Big( \inf_{s\in[0,T]} y^n_{s}(x) \le S \Big) \le \Prob\Big( \inf_{s\in[0,T]}\eta^n_s(x) \le S \Big).
\]
Here
\begin{align*}
\eta^n_t(x) &\ge \left( \e^{-C_a T} x - C_a T \right) + \wt\eta^n_t(x),
\end{align*} 
where 
\begin{align*}
\wt\eta^n_t(x) &= w^n_t(x) -w^n_{t^n_k}(x) + \sum_{j=\ov{0,k-1}} \e^{-C_a (t^n_{k+1} - t^n_{j+1})} \zeta^n_j, \\
&\quad t\in I^n_k, k =\ov{0,N^n-1},
\end{align*}
is a centered Gaussian process with 
\begin{align*}
\sup_{t\in [0,T]} \Var\left(\wt\eta^n_t(x)\right) \le  T.
\end{align*} 
The concentration inequality for a supremum of a Gaussian process~\cite[Theorem 2.1.1]{AdlTay09random} implies
\begin{equation}
\label{eq:measure.add.1}
\Prob\Big( \inf_{s\in[0,T]} y^n_{s}(x) \le S \Big) \le \exp\Big\{ -\frac{1}{2T} \Big( \e^{-C_a T}x - S - C_aT - \E \sup_{t\in [0,T]} \wt\eta^n_t(x) \Big)^2\Big\}. 
\end{equation}
Due to Dudley's entropy bound~\cite[Theorem 1.3.3]{AdlTay09random} for some absolute $K$
\begin{equation}
\label{eq:measure.add.2}
\E \sup_{t\in [0,T]} \wt\eta^n_t(x) \le K \int_0^\infty \left( \log M^n_\ve \right)^{\frac{1}{2}} d\ve,
\end{equation}
where $M^n_\ve$ is the smallest number of balls of size $\ve$  that cover  $[0,T]$ in the intrinsic metric 
\[
\rho^n(s_1, s_2) = \left(\E \left(\wt\eta^n_{s_1}(0) - \wt\eta^n_{s_2}(0)\right)^2 \right)^{\frac{1}{2}}
\]
 of the process $\wt\eta^n(0).$ Note that $\wt\eta^n(0)$ is a Wiener process on every $I^n_k$ and for $s_1 \in I^n_{k_1}, s_2 \in I^n_{k_2}, k_1 < k_2$ we have
\begin{align*}
\rho^n(s_1,s_2)^2 &= \sum_{j=\ov{0,k_1-1}} \e^{-2C_at^n_{j+1}} \left(\e^{-C_a t^n_{k_1+1}} - \e^{-C_a t^n_{k_2+1}} \right)^2 \delta^n_j \\
&+ \sum_{j=\ov{k_1+1,k_2-1}} \e^{-2C_a (t^n_{k_2+1} - t^n_{j+1})} \delta^n_j \\
&+ \left(1 - \e^{-C_a (t^n_{k_2+1}-t^n_{k_1+1})} \right)^2 (s_1 - t_{k_1}) \\
&+ \e^{-2C_a (t^n_{k_2+1}-t^n_{k_1+1})} (t^n_{k_1+1} -s_1) + s_2 - t_{k_2},
\end{align*}
so for a universal constant $C>1$
\begin{align*}
\rho^n(s_1,s_2) \le C\Big( (s_2-s_1)^{\frac{1}{2}} + \delta^n_{k_2} \left(  t^n_{k_1} \right)^{\frac{1}{2}}\Big).
\end{align*}
We assume $T=1$ for the rest of the proof. The diameter of $[0,T]$ in $\rho^n$ does not exceed $2C$ so $M^n_\ve = 1$ for $\ve\ge C$ for all $n.$ 

Assume $\ve > 2C (\delta_n)^{1/2}.$  For a unit interval, consider a $\frac{\ve^2}{4C^2}-$net $A$  in the Euclidean metric. Then for any $s$
\begin{align*}
\rho^n(s, A) &\le C \left( \frac{\ve}{2C} + (\delta_n)^{\frac{1}{2}}  \right) \le \ve, 
\end{align*}
so $A$ is an $\ve-$net for $[0,1]$ in the metric $\rho^n$ and 
\[
M^n_\ve \le \frac{4C}{\ve^2} + 1.
\]

Now assume that $\ve \le 2C (\delta_n)^{1/2}.$ We call an interval $I^n_k$ large, if $\ve \le 2C (\delta^n_k)^{1/2},$ and small, otherwise. We construct an $\ve-$net as follows. Let $B$ be a $\frac{\ve^2}{C_1}-$net for $[0,1]$ in the Euclidean metric where
\[
C_1 = 8(C^2 + 1).
\]
For each large $I^n_k$ let $B_k$ be a net of the same size such that $t^n_k\in B_k,$ also in the Euclidean metric. Set
\[
A_1 = B \cup \bigcup_{k=\ov{0,N^n-1}: I^n_k \mathrm{\ is\ large}} B_k. 
\]
 Since for large intervals $4C^2 \frac{\delta^n_k}{\ve^2} \ge 1,$ we have then
\begin{align*}
\sharp A_1 &\le \frac{C_1}{\ve^2} + 1 + \sum_{k=\ov{0,N^n-1}: I^n_k \mathrm{\ is\ large}} \Big(  \frac{C_1\delta^n_k}{\ve^2} + 1\Big) \\
&\le \frac{2C_1 + 4C^2}{\ve^2} + 1.
\end{align*}

We will show that $A_1$ is an $\ve-$net for $[0,1]$ in the metric $\rho^n.$ If $t \in I^n_k$ and $I^n_k$ is large, $\rho^n(t,A_1) \le \ve$ immediately. Assume $t\in I^n_k$ and $I^n_k$ is small. Then there exist numbers $0 \le m_1 \le k \le  m_2 \le N^n-1$ such that 
\begin{enumerate}
\item the intervals $I^n_{m_1}, \ldots, I^n_{m_2}$ are small;
\item either $m_1 =0$ or $I^n_{m_1-1}$ is large;
\item either $m_2 = N^n-1$ or $I^n_{m_2+1}$ is large;
\item at least one of $I^n_{m_1-1}$ and $I^n_{m_2+1}$ is large.
\end{enumerate}
Consider the case 
\[
L = t^n_{m_2+1} - t^n_{m_1} > \frac{2\ve^2}{C_1}.
\]
Then there exists $s\in B \cap \bigcup_{j = \ov{m_1, m_2}} I^n_j$ such that $|t-s| \le \frac{\ve^2}{C_1}$ and for some $j, m_1 \le j\le m_2 $
\[
\rho^n(s,t) \le C\Big( \frac{\ve}{C^{1/2}_1} + \left(\delta^n_{\max\{k,j\}} \right)^{1/2} \Big) \le  \ve.
\] 

If $L\le \frac{2\ve^2}{C_1},$ we consider two possibilities. If $I^n_{m_1-1}$ is large, then 
\begin{align*}
\rho^n(t, B_{m_1-1}) &\le \lim_{s\to t^n_{m_1}}\rho^n(t, s) + \frac{\ve}{C_1^{1/2}} \le C\left( L^{\frac{1}{2}} + \delta^n_k \right) + \frac{\ve}{C_1^{1/2}} \le \ve. 
\end{align*}
If $m_1=0,$ then $t^n_{m_2+1} = L,$ so
\begin{align*}
\rho^n(t, B_{m_2+1}) = \rho^n(t, t^n_{m_2+1}) \le C \left( (t^n_{m_2+1} - t)^{\frac{1}{2}} + L^{\frac{1}{2}}\right) \le 2C L^{\frac{1}{2}} \le \ve. 
\end{align*}
This proves the claim.

Estimating the integral in~\eqref{eq:measure.add.2} we obtain for some $K_1$
\[
\sup_{n\in\mbN} \E \sup_{t\in [0,T]} \wt\eta^n_t(x)  \le K \int_0^{C} \Big( \log\Big( \frac{K_1}{\ve^2} + 1\Big) \Big)^{\frac{1}{2}}d\ve  < \infty,
\]
so \eqref{eq:measure.add.1} yields
\[
 \sup_{x \ge x_0} \sup_{n\in \mbN}\Prob\Big( \inf_{s\in[0,T]}\eta^n_s(x) \le S \Big) \le C_0 \e^{-\frac{x^2}{C_0}}.
\]
 for  sufficiently large absolute  $x_0$ and  $C_0.$  Since the same estimate can be obtained for  $\Prob( \sup_{t\in[0,T]}\eta^n_t(x) \ge S)$ for  $n\in \mbN$ and  negative $x,$   \eqref{eq:add.22} follows.

For the limit process,  we get for $x \gg S$
\[
\Prob\Big( \inf_{s\in[0,T]} X^{\vf,a}_{0,s}(x) \le S \Big) \le \Prob\Big( \inf_{s\in[0,T]}\xi_s(x) \le S \Big), 
\]
where 
\[
d\xi_t(x) = - C_a(1 + \xi_t(x)) dt + dw_t, \quad \xi_0(x) =x,
\]
$w$ being a Wiener process. Since
\[
\xi_t(x) =  x\e^{-C_at} + \e^{-C_at} - 1 +\e^{-C_at} \int_0^t \e^{C_a s} dw_s, \quad t \ge 0, 
\]
where the last term is a continuous square integrable martingale with bounded quadratic variation,  and  $\Prob( X^{\vf,a}_{0,t}(x) \ge -S)$ can be estimated similarly, \eqref{eq:add.21} follows. 
\qed

{\it Proof of Corollary \ref{cor:add.2}.} By Theorem  \ref{th:add.1} $\mu^n_t \Rightarrow \mu_t, n\to\infty,$ vaguely. Since all measures are probabilistic, the convergence also holds in the weak topology~\cite[Theorem 4.9]{Kal83Random}. \qed

\section{Convergence of dual flows}
We assume that $((u^n,y^n))_{n\in\mbN}=\mathrm{Spl}(X^{\vf,0}; a; \mathcal{T})$ and establish the convergence of so-called dual flows~\cite{Ha84Coalescing, Ria18Duality, Vov18Convergence, AmaTaYu19Convergence}. 
Here given a coalescing stochastic flow $X$ the dual flow $\wh X$ is defined via
\begin{align}
\label{eq:dual.def}
\wh X_{s,t}(x) &= \inf \left\{y \in \mbR\mid X_{T-t,T-s}(y) >x \right\} \notag \\
&=  \inf\big\{X_{r,T-t}(y)\mid X_{r,T-s}(y)> x, y\in\mbR, r \in[0,T-t]\big\}, \notag  \\
& s,t\in[0,T], s\le t, \ x \in \mbR,
\end{align}  
and is again a collection of $D^\uparrow(\mbR)-$valued random elements.

To start, we need one extension of the splitting scheme~\eqref{eq:splitting.flows}. 
 Let $l(s)$ equal the unique $k$ such that $s\in I^n_{k}, s\in [0,T].$
Define
\begin{align*}
t^{n}_{k}(s) &= \max\{s, t^{n}_{k}\}, \quad k =\ov{l(s),N^n-1}, n\in\mbN, \notag \\ 
u^n_{s,t}(x) &= F_{t - t^{n}_{k}(s)}\left(y^n_{s,t^{n}_{k}-}(x)\right), \notag \\
y^n_{s,t}(x) &= X^{\vf,0}_{t^{n}_{k}(s),t}\left(u^n_{s,t^{n}_{k+1} -}(x)\right), \notag \\
y^n_{s,s-}(x) &= x, \\
& \quad t \in I^n_{k}, k =\ov{l(s),N^n-1}.
\end{align*} 

For any $s\in[0,T)$ and any $f\in D([s,T],\mbR)$ we extend $f$ onto $[0,T]$ by setting 
\[
f^e_r = f_r\1\left[r \in [s,T]\right] + f_s\1\left[r \in [0,s)\right].
\]
For instance,  $y^{n,e}_{s,\cdot}(x), u^{n,e}_{s,\cdot}(x)$ are random elements in $D([0,T],\mbR).$

Since constructing the families $\{y^n_{s,\cdot}(x)\mid x\in\mbR\}$ uses the single flow $X^{\vf,0}$ for all $s,$ the mappings $\{y^n_{s,t}\mid 0\le s\le t \le T\}$ are consistent  and form a coalescing flow, so using~\eqref{eq:dual.def}, one defines   the corresponding dual flow $\wh{y}^n=\{\wh y^n_{s,t} \mid 0\le s\le  t \le T\}.$

We use the idea from \cite{Vov18Convergence}, where a precise construction of the dual flow as a function that preserves the weak convergence is given.  

Consider the set $\{(s_n,x_n) \in [0,T]\times\mbR \mid n\in\mbN\}$ containing all points  whose coordinates are dyadic numbers. The corresponding version of Theorem~\ref{th:con.1} implies that for any $m\in\mbN$
\[
Y^{n,m} = \left( y^{n,e}_{s_1,\cdot}(x_1), \ldots, y^{n,e}_{s_m,\cdot}(x_m)\right) \Rightarrow X^m=\left( X^{\vf,a,e}_{s_1,\cdot}(x_1), \ldots,  X^{\vf,a,e}_{s_m,\cdot}(x_m) \right), \quad n\to\infty,
\]
in $D([0,T],\mbR^m).$

Denote by $P_k$ the projector on the first $k$ coordinates in the space $D([0,T],\mbR)^\infty$ endowed with the product topology and by $Q_k$ the projector on the $k-$th coordinate in the same space. Consider the  $D([0,T],\mbR)^\infty-$valued random elements $X, \wh X, Y^n, \wh Y^n, n\in\mbN,$ defined via
\begin{align*} 
P_m(Y^n) &= Y^{n,m}, \\
 P_m(X) &=X^m, \\
P_m(\wh Y^n) &= \left(\wh y^{n,e}_{s_1,\cdot}(x_1),\ldots, \wh y^{n,e}_{s_m,\cdot}(x_m) \right), \\
 P_m(\wh X) &= \Big(\wh X^{\vf,a,e}_{s_1,\cdot}(x_1),\ldots, \wh X^{\vf,a,e}_{s_m,\cdot}(x_m) \Big), \\ 
& n, m\in\mbN.
\end{align*}

As in \cite[p. 86]{Vov18Convergence}, we obtain 
\[
Y^n \Rightarrow X, \quad n\to\infty,
\]
in $D([0,T],\mbR)^\infty.$

Consider a mapping $I\colon D([0,T],\mbR)^\infty \mapsto D([0,T],\mbR)^\infty$ defined via
\begin{align*}
Q_j I(\psi)_r &= \inf\{ Q_i \psi_r \mid Q_i \psi_T\ge x_j, s_i \le r \}, \quad r\in[s_j,T], \\
Q_j I(\psi)_r &=  Q_j I(\psi)(s_j),  \quad  r\in[0,s_j),  \\
& j\in\mbN,  \ \psi \in D([0,T]^\infty).
\end{align*}
Then one can check that  $I(Y^n)=\wh Y^n$ a.s., $n\in\mbN,$ and $I(X)=\wh X$ a.s..

\begin{definition}
Let $D_1$ be a set of $\psi\in D([0,T],\mbR)^\infty$ such that 
\begin{enumerate}
\item $Q_k{\psi}_{s_k} = x_k, \ k\in\mbN;$
\item if for some $ j_1, j_2\in\mbN$ and some $s\in[0,T]$ $Q_{j_1}{\psi}_s\ge Q_{j_2} {\psi}_s,$ then $ Q_{j_1} {\psi}_t\ge Q_{j_2} {\psi}_t$ for  $t\in [s;T].$
\end{enumerate} 
\end{definition}

\begin{definition}
Let $D_2$ be a subset of $D_1$ such that for any $ \psi\in D_2$
\begin{enumerate}
\item  $\forall k\in\mbN \ \exists \ve_k>0$  
\begin{align*}
\label{prop:divergence}
& \forall i\in\mbN\colon \big(x_i \ge Q_k I(\psi)_{s_i}\big) \Rightarrow \big(Q_i \psi_{T} - x_{k} \ge \ve_{k}\big), \\
& \forall i\in\mbN\colon \big(x_i < Q_k I(\psi)_{s_i}\big) \Rightarrow \big(x_{k} - Q_i \psi_T \ge \ve_{k}\big);
\end{align*}
\item 
$\forall \delta>0 \ \forall M>0 \ \exists L\in\mbN$ 
\begin{align*}
\sup_{l=\overline{0,\lceil T\rceil 2^{L}}}\sup_{j\colon |x_j|\le M, s_j = l2^{-L}} \sup_{\tau\in[0,2^{-L}]} \big| Q_j\psi_{s_j+\tau}- x_j \big| \le \delta.
\end{align*}
\end{enumerate}
\end{definition}

\begin{lemma}[{\cite[pp. 86--87]{Vov18Convergence}}]
\label{lem:dual.1}
Assume that $\psi_n\rightarrow\psi, n\rightarrow\infty,$ in $D([0,T],\mbR)^\infty,$ $\psi_n \in D_1, n\in\mbN, \psi\in D_2,$ and $Q_j \psi \in C([0,T]), j\in\mbN.$
Then $I(\psi_n)\rightarrow I(\psi), n\rightarrow\infty,$ in $D([0,T],\mbR)^\infty.$ 
\end{lemma}

\begin{lemma}
\label{lem:dual.2}
Assume that $\vf\in \Phi_\alpha$ and $a \in A_\beta$ with $\beta - \alpha > -  1, \alpha< 2.$ Then $Y^n \in D_1$ a.s., $n\in\mbN.$ $X \in D_2$ a.s.. 
\end{lemma}
\proof Proceeding  as in~{\cite[p. 87]{Vov18Convergence}}, one uses~\eqref{eq:app.1} in the proof of Proposition~\ref{prop:4.7} and the fact that the set $\{X^{\vf,a}_{s,t}(x) \mid x\in\mbR\}$ is a.s. locally finite by Theorem~\ref{th:flow.coalescence} for any $s,t, s<t.$ \qed

Combining Lemmas \ref{lem:dual.1} and \ref{lem:dual.2} yields the following result.

\begin{theorem}
\label{th:dual.1}
Assume that $\vf\in \Phi_\alpha$ and $a \in A_\beta$ with $\beta - \alpha > -  1, \alpha< 2.$ Then $\wh Y^n \Rightarrow \wh X, n\to\infty,$ in $D([0,T],\mbR)^\infty.$ In particular, for any $t_1, \ldots, t_m,$ $v_1, \ldots, v_m$ and $m\in\mbN$ 
\[
\left(\wh y^{n,e}_{t_1,\cdot}(v_1), \ldots, \wh y^{n,e}_{t_m,\cdot}(v_m) \right) \Rightarrow \left( \wh X^{\vf,a,e}_{t_1,\cdot}(v_1), \ldots,  \wh X^{\vf,a,e}_{t_m,\cdot}(v_m) \right), \quad n\to\infty, 
\]
in $D([0,T],\mbR^m).$
\end{theorem}

\section{Estimates for flows  with  Holder continuous $\sqrt{1-\vf}$}
$((u^n,y^n))_{n\in\mbN}=\wt{\mathrm{Spl}}(X^{\vf,a}; a; \mathcal{T})$ and $W=\cW(X^{\vf, a})$ are considered in this section.

Assume that $\vf\in C^2(\mbR),$ which is equivalent to the finiteness of $\vf^{\prime\prime}(0).$ As noted in~\cite[\S 3]{Ha84Coalescing}, $\sqrt{1-\vf}\in Lip(\mbR)$ then. Thus results of~\cite{BenGloRas92Approximations} are applicable. $X^{\vf,a}$ is a flow of homeomorphisms by \cite[Theorem 4.5.1]{Ku90Stochastic}. 

\begin{theorem}[{\cite[Corollary 4.2]{BenGloRas92Approximations}}]
\label{th:homeo.1}
For any $M \ge 0$ and some $C = C(M, T)>0$
\begin{align*}
\sup_{x\in [-M,M]} &\E \sup_{t\in [0,T]} \left( y^n_t(x) - X^{\vf,a}_{0,t}(x)\right)^2 \le C \delta_n, \\
\sup_{x\in [-M,M]} &\sup_{t\in [0,T]} \E\left( u^n_t(x) - X^{\vf,a}_{0,t}(x)\right)^2 \le C \delta_n.
\end{align*} 
\end{theorem}


It is possible to refine the order of convergence for $u^n.$
\begin{proposition}
\label{prop:homeo.1}
For any $M \ge 0$ and some $C = C(M, T)>0$ 
\begin{align*}
\sup_{x\in [-M,M]} \E \sup_{t\in [0,T]} \left( u^n_t(x) - X^{\vf,a}_{0,t}(x)\right)^2 \le C \delta_n \log \delta_n^{-1},
\end{align*}
\end{proposition}
\proof  
Dropping the $x$ argument and setting  
\[
m_t = \E \sup_{s \in[0,T]}\max \Big\{ \left(X^{\vf,a}_{0,s} - u^n_s\right)^2, \left(X^{\vf,a}_{0,s} - y^n_s\right)^2  \Big\}.
\]
we get by Lemma \ref{lem:hol.3} for some $C_1,C_2$
\begin{align*}
\label{eq:homeo.1}
m_t &\le C_1\Big( \int_0^{t} m_s ds + \E \sup_{s \in[0,T] } \Big( \left( r^n_s \right)^2 + \left( l^n_s \right)^2\Big) \Big)\\
 &\le C_2\Big( \int_0^{t} m_s ds + \delta_n \log \delta_n^{-1} \Big).
\end{align*}
\qed 

The order of convergence in Proposition \ref{prop:homeo.1} cannot be improved, as shown by the following example.

\begin{example}
Let $a\equiv 0, T = 1$ and $t^{n}_{k} = \frac{k}{n}.$ Then $u^n_t(x)=X^{\vf,0}_{0,\frac{k}{n}}(x)$ on $I^n_{k},$ so
\begin{align*}
\E \sup_{k=\ov{0,n-1}} \sup_{t\in [\frac{k}{n}, \frac{k+1}{n})} \Big(X^{\vf,0}_{0,t}(x) -u^n_t(x)\Big)^2 &= \frac{2}{n} \E \sup_{k=\ov{1,n}} \eta_k^2, 
\end{align*}
where $\eta_k, k=\ov{1,n},$ are independent $\cN(0,1)$ random variables. It is well known that
\[
\lim_{n\to \infty} \frac{\E \sup_{k=\ov{1,n}} \eta_k}{(2 \log n)^{1/2}} =1, 
\]
 so
\[
\liminf_{n\to\infty} \frac{\E  \sup_{t\in [0,T]} \Big(X^{\vf,0}_{0,t}(x) -u^n_t(x)\Big)^2}{4n^{-1} \log n} \ge 1.
\]
\end{example}

The following  fact is well known (e.g.~\cite{DoFo16Rate, DorVov20ApproximationsEng}). 

\begin{proposition}
\label{prop:inv.mapping.measure}
Let $f$ be \cadlag and non-decreasing on $[0,1].$ Set $\vk = \Leb\vert_{[0,1]}\circ f^{-1}$ and $F(\cdot) = \vk([0,\cdot]).$ Then $F^{-1} = f$ on $[0,1].$
\end{proposition}

\begin{theorem}
For some $C = C(T)>0$
\[
\sup_{t\in[0,T]} W_{1,2}\left( L_t, L^n_t \right) \le C \delta_n, \quad n\in\mbN.
\]
\label{th:homeo.2}
\end{theorem}
\proof By \cite[Remark 2.19]{VillaniTopics03} 
\[
W_{1,2}\left( L_t, L^n_t \right) = \E \int_0^1 \left( F_t^{-1}(u) -F_{n,t}^{-1}(u)  \right)^2 du,
\]
where $F_t^{-1}, F_{n,t}^{-1}$ are generalized \cadlag inverses of 
\begin{align*}
F_t(u) &= \mu_t((-\infty,x]) = \Leb\left\{ u\in [0,1] \mid X^{\vf,a}_{0,t}(u) \le x\right\}, \\
F_{n,t}(u) &= \mu^n_t((-\infty,x]) = \Leb\left\{ u\in [0,1] \mid y^n_{t}(u) \le x\right\},
\end{align*} 
respectively. Thus, by Proposition \ref{prop:inv.mapping.measure},
\[
W_{1,2}\left( L_t, L^n_t \right) = \E \int_0^1 \left( X^{\vf,a}_{0,t}(u) -y^n_t(u)  \right)^2 du,
\] 
and the application of Theorem \ref{th:homeo.1} establishes the claim. \qed

Now assume that $\sqrt{1-\vf} \in H_\beta(\mbR)$ for some $\beta \in[\frac{1}{2},1)$.  The Euler-Maruyama scheme for an SDE with such coefficients was considered in \cite{GyRa11Note}, where a suitable  modification of the Yamada-Watanabe method was developed. Essentially, recalling \eqref{eq:notation.1}, using Lemmas~\ref{lem:hol.2} and~\ref{lem:hol.3} to estimate $\E |r^n_s(x)|^\beta$ and $\E |l^n_s(x)|^\beta,$  one proceeds by repeating the reasoning for \cite[Proposition 2.2]{GyRa11Note} line by line to draw the following two conclusions.

\begin{theorem}
\label{th:hol.2}
For some $C=C(x,\beta, T)>0$
\begin{align*}
\E \sup_{t\in [0,T]}\left(X^{\vf,a}_{0,t}(x) - y^n_t(x)\right)^2 &\le \frac{C}{\log \delta_n^{-1}}, \quad \beta = \frac{1}{2}, \\
\E \sup_{t\in[0,T]} \left(X^{\vf,a}_{0,t}(x) - y^n_t(x)\right)^2 &\le C \delta_n^{\beta-\frac{1}{2}}, \quad \beta \in \Big(\frac{1}{2},1\Big). 
\end{align*}
\end{theorem}
\begin{theorem}
\label{th:hol.3}
For some $C=C(\beta, T)>0$
\begin{align*}
\sup_{t\in[0,T]} W_{1,2}\left( L_t, L^n_t \right)  &\le \frac{C}{\log \delta_n^{-1}}, \quad \beta = \frac{1}{2}, \\
\sup_{t\in[0,T]} W_{1,2}\left( L_t, L^n_t \right) &\le C \delta_n^{\beta-\frac{1}{2}}, \quad \beta \in \Big(\frac{1}{2},1\Big).
\end{align*}
\end{theorem}

\begin{remark}
Following assumptions in~\cite{GyRa11Note}, consider $a =a_1 + a_2, a_1\in Lip(\mbR), a_2\in H_{\alpha}(\mbR)$ for some $\alpha\in (0,1)$ and assume that $a_2$ is non-increasing. Theorems~\ref{th:hol.2} and~\ref{th:hol.3} can be extended to such $a$ as follows. Consider the extension of the splitting scheme given in Remark~\ref{rem:weak.1} with $\ve_n = \tfrac{1}{n^{1/3}\log \delta_n^{-1}}$ for $\beta =\tfrac{1}{2}$ and $\ve_n = \frac{1}{n^{\frac{1}{2}+\min\{\frac{\alpha}{2}, \beta-\frac{1}{2}\}}}$ for $\beta\in(\tfrac{1}{2},1).$ Then
\begin{align*}
\E \sup_{t\in [0,T]}\left(X^{\vf,a}_{0,t}(x) - y^n_t(x)\right)^2 &\le \frac{C}{\log \delta_n^{-1}}, \quad \beta = \frac{1}{2}, \\
\E \sup_{t\in[0,T]} \left(X^{\vf,a}_{0,t}(x) - y^n_t(x)\right)^2 &\le C \delta_n^{\min\{\frac{\alpha}{2}, \beta-\frac{1}{2}\}}, \quad \beta \in \Big(\frac{1}{2},1\Big). 
\end{align*}
\end{remark}

\begin{appendices}
\section{Existence of Harris flows with drift}

Assume that $\vf\in\Phi_*$ and a measurable $a$ with linear growth are fixed.


Consider the following operator acting on $C^2(\mbR^n)$
\begin{equation*}
\label{eq:main.operator}
\cA_n = \frac{1}{2}\sum_{k,j=\ov{1,n}} \vf(x_k-x_j) \frac{\pt^2}{\pt x_k \pt x_j} + \sum_{k=\ov{1,n}} a(x_k) \frac{\pt }{\pt x_k}. 
\end{equation*}
The degeneracy of the matrix $\|\vf(x_k-x_j)\|_{k,j=\ov{1,n}}$ on the boundary of $D_n = \{x\mid  x_1 < \ldots < x_n\}$  requires an extension of the results presented in \cite{Stroock1Va97Multi}, which was directly discussed by Harris.  However, \cite[\S 1.9-13]{Pin95Positive} later provided a rigorous framework in terms of generalized martingale problems in domains, and we adopt this approach.

Consider $C(\mbR_+, \mbR^n)$ endowed  with the Borel $\sigma-$algebra and the topology of uniform convergence on compact sets. Define $\mcC_n$ be the set of elements of $C(\mbR_+, \mbR^n)$ whose coordinates merge after a meeting. 

\begin{definition}(\cite[Definition 2.1]{Ha84Coalescing})
\label{defn:C.solution}
A family $P_x, x \in \mbR^n,$ is called a $\mcC_n-$solution if it solves the martingale problem for $\cA_n$ in $\mbR^n$ (in the sense of~\cite{Stroock1Va97Multi}) and for each $x$ $P_x(\mcC_n) = 1$ and $P_x(\omega\mid \omega_0 =x\})=1.$
\end{definition}

The solution of the generalized martingale problem for  $\cA_n$ does not explode to infinity, so 
one can use solutions of generalized martingale problems to construct $\mcC_n-$solutions in \cite[Lemma 3.2]{Ha84Coalescing}. The next proposition is essentially \cite[Lemmas 2.2, 3.2]{Ha84Coalescing}, with the exception of the conclusion about the measurability of the mapping $x \mapsto P_x,$ which follows from  \cite[Proof of Theorem 4.4.6]{EthKur09Markov}, and the conclusion about the Feller property, which follows from \cite[Theorem 1.13.1]{Pin95Positive}. 

\begin{proposition}
\label{3.}
There exists a unique $\mcC_n-$solution $P_x, x\in \mbR^n.$  This solution is strong Markov,  Feller and such that the mapping $x \mapsto P_x \in \cM_0(C(\mbR_+, \mbR^n))$ is measurable, where $\cM_0(C(\mbR_+, \mbR^n))$ is the space of probability measures on $C(\mbR_+, \mbR^n)$ endowed with the topology of weak convergence.
\end{proposition}

To construct $X^{\vf,a}$ in Theorem~\ref{th:flow.existence}, one can  repeat the reasoning in~\cite[\S 4]{Ha84Coalescing}, the only missing ingredient that has to be checked directly being the following extension of \cite[Theorem 4.7]{Ha84Coalescing}.  

\begin{proposition}
\label{prop:4.7}
For any $t>0$ the family $\{X^{\vf, a}_{0,s}\mid s \in [0,T]\}$ is a.s. uniformly continuous in $D^{\uparrow}(\mbR)$ endowed with the topology of uniform convergence on compact sets.
\end{proposition} 
The next lemma is a straightforward corollary of the Gronwall lemma. 
\begin{lemma}
\label{lem:app.1}
For $C$ such that  $|a(x)| \le C(1 + |x|), x\in \mbR,$  and a standard Wiener process~$w$ 
\begin{align*}
\E \sup_{s\in[0,T]}\left| X^{\vf,a}_{0,t}(x) - x\right|^{p} \le 2^{p-1} \e^{pCT} \Big( C^{p} (1 + |x|)^{p} t^{p} + \E \sup_{s\in[0,T]} |w(s)|^{p} \Big), \quad p \ge 1.
\end{align*}
\end{lemma}

{\it Proof of Proposition \ref{prop:4.7}.}
We modify the original proof. Define, for fixed $T>1$ and $N > 1$ 
\[
U_{knN} = \sup_{t\in [k2^{-n}T,(k+1)2^{-n}T]} \sup_{|x| \le N} \left| X^{\vf,a}_{0,t}(x) - X^{\vf,a}_{0,k2^{-n}T}(x) \right|, \quad k,n\in\mbN.
\]
It is sufficient to show for any $\ve > 0$
\begin{equation}
\label{eq:app.1}
\sum_{n\ge 1} \Prob\Big( \bigcup_{k=\ov{0,2^n-1}} \{U_{knN} \ge \ve\} \Big) \le \sum_{n\ge 1} \sum_{k=\ov{0,2^n-1}} \Prob\left( U_{knN} \ge \ve \right)< \infty.
\end{equation}
Proceeding as in \cite[Proof of Lemma 4.5]{Ha84Coalescing} and using Lemma \ref{lem:app.1}, we have, for any $m\ge 1,$ $a,b\in\mbR, b-a>1, $ and some $C> 0$
\begin{align*}
&\Prob\Big( \sup_{s\in[0,T], x\in[a,b]} \left|X^{\vf,a}_{0,s}(x) - x\right|\ge \ve \Big) \\
&\qquad\qquad \le 2^{m+1}(b-a) \sup_{x\in[a,b]} \Prob\Big( \sup_{s\in[0,T]}\left|X^{\vf,a}_{0,s}(x) - x\right|\ge \ve -2^{-m}\Big) \\
&\qquad\qquad \le \frac{2^{m+1}(b-a)}{(\ve - 2^{-m})^4} \sup_{x\in[a,b]} \E \sup_{s\in[0,T]}\left(X^{\vf,a}_{0,s}(x) - x\right)^4 \\
&\qquad\qquad \le  C 2^m  \frac{(b-a)^5 t^4 + (b-a) t^2}{(\ve - 2^{-m})^4}.
\end{align*}
Thus, for any fixed $m$ such that $2^{-m} \le \tfrac{\ve}{2}$ and some $C_1 = C_1(m,N,T)>0,$
\begin{align*}
 \Prob\left( U_{knN} \ge \ve \right) &\le \E\Prob \Big( \sup_{t \in[0,2^{-n}T]}\sup_{y \in A} \left| X^{\vf,a}_{0,t}(x) - x\right| \ge \ve  \Big)\Big|_{A = \left[X^{\vf,a}_{0, k2^{-n}T}(-N),X^{\vf,a}_{0, k2^{-n}T}(N)\right]} \\ 
 &\le C \ve^{-4} 2^{-2n},
\end{align*}
which \eqref{eq:app.1} follows from.
\qed

This concludes the proof of Theorem \ref{th:flow.existence}.

Assume $\vf\in\Phi_{\alpha},$ $a \in A_\beta,$  $\beta -\alpha> -1$ and $\alpha<2.$ Theorem \ref{th:flow.coalescence} corresponds to \cite[Theorem 7.4]{Ha84Coalescing}, the crucial point being the estimate %
\begin{align}
\label{eq:app.2}
\E \sharp \left\{ X^{\vf, 0}_{0,t}(u) \mid u \in [0,M] \right\} &\le 1 + \limsup_{n\to\infty}\sum_{k=\ov{0,n-1}} \Prob\Big( X^{\vf,0}_{0,t}\Big(\frac{k+1}{n}M\Big) >X^{\vf,0}_{0,t}\Big(\frac{k}{n}M\Big) \Big), \notag  \\
& \quad M \in\mbR_+, n\in\mbN,
\end{align}    
where for some $C=C(t)$ and $0<y-x$ 
\begin{align}
\label{eq:app.3}
\Prob\Big( X^{\vf,0}_{0,t}(y) >X^{\vf,0}_{0,t}(x) \Big) \le C (y-x),
\end{align}
so then 
\[
\E \sharp \left\{ X^{\vf, 0}_{0,t}(u) \mid u \in [0,M] \right\}  \le 1 + CM. 
\]
Since \eqref{eq:app.2} holds for $X^{\vf,a},$ too, it is sufficient to prove \eqref{eq:app.3} in the case of non-zero drift. 

In what follows we refer to \cite[\S 5.5]{KaShre91Brownian} and \cite[Chapter 16]{Brei68Probability} for the theory of Feller's one-dimensional diffusions and a classification of boundaries and to \cite{GoJaYor03Bessel} for basic facts about Bessel processes including those with negative dimension.    

Let $\xi$ be a diffusion on $\mbR_+$ with generator 
\[
(1-\vf(x))\frac{d^2}{dx^2} + \rho(x)\frac{d}{dx}
\]
and an absorbing boundary at $0.$ 
The scale function can be defined as 
\[
p(x) = \int_0^x \e^{-\int_0^y \frac{\rho(z)}{1-\vf(z)}dz} dy,
\]
so the speed measure is
\[
m(dy) = \frac{\e^{\int_0^y \frac{\rho(z)}{1-\vf(z)} dz}}{1-\vf(y)} dy.
\]
Set $\wt{C} = \min\{\wt{C}_\rho, \wt{C}_\vf\}$ and define
\[
v(x) = \int_{\wt{C}}^{x} (p(x) - p(y)) m(dy) \ge 0, \quad x\in\mbR_+.
\]
Since  
\[
\lim_{\ve\to 0+}v(\ve) \le C_\vf^{-1} \e^{\int_0^{\wt{C}} \frac{\rho(y)}{1-\vf(y)}dy}  \limsup_{\ve\to 0+} \int_\ve^{\wt{C}} \frac{y-\ve}{y^\alpha} dy <\infty,
\]
the boundary $0$ is accessible in finite time with positive probability $p_0.$ 
  It is worth noting that $p_0<1$ in general, contrary to the case $\rho\equiv 0,$ as shown below.

\begin{example}
Let $\rho(x)=0.$ Then $\lim_{x\to\infty} p(x)=\infty,$ $m([1,\infty))=\infty,$  so $\infty$ is a natural boundary and the diffusion $\xi$ hits $0$ a.s. in finite time $\tau,$ though $\E \tau = \infty.$ Indeed, $\xi = \xi_0 + {w}_{\langle \xi \rangle}$ for some Wiener process $w$, where $\langle \xi \rangle_t \le 2t, t\ge0,$ so, if $\tau_y = \inf\{s \mid \xi_s = y\}, y\in\mbR_+,$ we have   
\begin{align*}
\E \tau &> \Prob\left( \tau_{\xi_0+1} < \infty \right) \E_{\xi_0+1} \tau_{\xi_0}  \\
&=  \Prob\left( \tau_{\xi_0+1} < \infty \right) \int_0^{\infty} \Prob\Big( \xi_0+1 + \inf_{s\in[0, t]} w_{\langle \xi \rangle_s} > \xi_0 \Big) dt \\
&\ge  \Prob\left( \tau_{\xi_0+1} < \infty \right) \int_0^{\infty} \Prob\Big( \inf_{s\in[0,2t]} w_s > -1 \Big) dt \\
&= +\infty.
\end{align*}
In fact,  $0$ is an exit for $\alpha\in [1;2)$ and a regular boundary for $\alpha\in(0;1).$   
\end{example}

\begin{example}
\label{ex:app.2}
Let $\rho(x) =x.$ Since $-\tfrac{z}{1-\vf(z)} \le -z$ on $\mbR_+$, 
\[
\lim_{x\to\infty}p(x) \le \int_{0}^{\infty} \e^{- \frac{y^2}{2}} dy <\infty,
\]
so $\infty$ is accessible. Then
\[
p_0 \le  \frac{\lim_{x\to\infty} p(x) - p(\xi_0)}{\lim_{x\to\infty} p(x)} < 1. 
\]
\end{example}

\begin{example}
Let $\rho(x) =1$ and $\alpha<1.$ Then the situation of Example \ref{ex:app.2} repeats. In particular, given $y>x$ the difference $X^{\vf,a}_{0,\cdot}(y)-X^{\vf,a}_{0,\cdot}(x)$ goes to infinity with a positive probability in either example.
\end{example}

The next lemma is proved with a standard localization argument.
\begin{lemma}
\label{lem:app.2}
For any $(u_1,u_2) \in D_2$ and any $t >0$
\[
\Prob\left(X^{\vf,a}_{0,t}(u_2) > X^{\vf,a}_{0,t}(u_1) \right) \le \Prob_{u_2-u_1} \left( \xi_t > 0\right).
\]
\end{lemma}

\begin{remark}
By the localization argument, one can prove a stronger statement: given $\xi_0=u_2-u_1>0$ we have $\xi \ge X^{\vf,a}_{0,\cdot}(u_2) -X^{\vf,a}_{0,\cdot}(u_1)$ a.s. if $\xi$ solves, up to the moment of hitting $0,$ 
\[
d\xi_t = \rho(\xi_t) dt + \left(2(1-\vf(\xi_t))\right)^{1/2} dw_t,
\] 
where the Wiener process $w$ is such that
\begin{multline*}
X^{\vf,a}_{0,t}(u_2)-X^{\vf,a}_{0,t}(u_1)  \\ = \int_{0}^t \left[a\left(X^{\vf,a}_{0,s}(u_2)\right)-a\left(X^{\vf,a}_{0,s}(u_1)\right)\right]ds + 2^{1/2}\!\! \int_0^t \! \left[1-\vf\left(X^{\vf,a}_{0,s}(u_2)-X^{\vf,a}_{0,s}(u_1)\right)\right]^{1/2} dw_s.
\end{multline*}
\end{remark} 

\begin{proposition}
Suppose $\vf\in \Phi_\alpha,$ $a \in A_\beta,$ $\beta - \alpha > -  1$ and $\alpha<2.$ Then for any $t>0$ and any $x,y \in \mbR, y >x,$ 
\begin{align}
\label{eq:app.3.1}
\Prob\Big( X^{\vf,a}_{0,t}(y) >X^{\vf,a}_{0,t}(x) \Big) \le C (y-x),
\end{align}
for some $C=C(t).$

\end{proposition}
\proof By Lemma \ref{lem:app.2} it is sufficient to prove \eqref{eq:app.3.1} for the diffusion $\xi.$ We follow the idea of~\cite{Ha84Coalescing} of switching to a squared Bessel process. Since $\xi$ is not in the natural scale, one step is added and the coefficients are distorted, so we provide necessary details.  

Set $\sigma = (2(1-\vf))^{1/2}, p(\infty) = \lim_{x\to\infty} p(x).$  The diffusion $\wt\xi= p(\xi)$ on $(0,p(\infty))$ with absorption at $0$ has  generator $\frac{1}{2}\wt\sigma^2(x)\frac{d^2}{dx^2}$, where 
\begin{align*}
\wt\sigma(y) &=\sigma(p^{-1}(y)) p^\prime(p^{-1}(y))= \sigma(p^{-1}(y)) \e^{- \int_0^{p^{-1}(y)} \frac{\rho(z)dz}{1-\vf(z)}}, \\
p^{-1}(y) &= \int_0^y \frac{du}{p^\prime(p^{-1}(u))}= \int_0^y \e^{ \int_0^{p^{-1}(u)} \frac{\rho(z)dz}{1-\vf(z)} } du, \\
&y\in [0,p(\infty)).
\end{align*}
Set $\wt\theta =\inf\{s\mid \wt\xi_s =0\}.$ Define 
\begin{align*}
\psi(y) &= \frac{\wt\sigma^2(y)}{y^\alpha}, \quad y \in [0,p(\infty)), \\
\tau(t) &= \inf\Big\{s \mid \int_0^s \psi(\wt\xi_r) dr = t\Big\}, \quad t \le \int_0^{\wt\theta} \psi(\wt\xi_s)ds. 
\end{align*} 
Then $\eta = \wt\xi_{\tau}$ is a diffusion on $(0,p(\infty))$ with generator $\frac{1}{2}y^{\alpha} \frac{d^2}{dy^2}$ and absorption at $0,$ and $\wt\eta = (\frac{2}{2-\alpha})^2\eta^{2-\alpha}$ is a diffusion on $(0,(\frac{2}{2-\alpha})^2p(\infty)^{2-\alpha})$ with generator $2y\frac{d^2}{dy^2} + \delta \frac{d}{dy}$ and absorption at $0.$ Here $\delta = \frac{2(1-\alpha)}{2-\alpha} \in (-\infty,1).$ Note that the diffusion $\wt{\wt\eta}$  on $[0,\infty)$ with generator $2y\frac{d^2}{dy^2} + \delta \frac{d}{dy}$ is a squared Bessel process with  dimension $\delta$ and always hits $0.$ It is known~\cite[Example 6.4]{KaRuf16Distribution} that   
\begin{equation}
\label{eq:bessel.time}
\Prob_{\wt{\wt\eta}_0}\left( \wt{\wt\eta}_t > 0 \right) = \frac{1}{\Gamma(\frac{1}{2-\alpha})} \int_0^{{\wt{\wt\eta}_0}/(2t)}  s^{-\delta/ 2} \e^{-s} ds \le \frac{2-\alpha}{\Gamma(\frac{1}{2-\alpha})} \Big( \frac{\wt{\wt\eta}_0}{2t} \Big)^{\frac{1}{2-\alpha}}.
\end{equation}
  
If $p(\infty) = \infty$ we have, for any $\eta_0>0,$ 
\begin{equation}
\label{eq:app.bessel.est}
\Prob_{\eta_0}\left( \eta_t > 0 \right) = \Prob_{(\frac{2}{2-\alpha})^2\eta_0^{2-\alpha}}\left( \wt{\wt\eta}_t > 0 \right),
\end{equation} 
but if $p(\infty) < \infty$
\[
\Big\{ {\wt\eta} \mathrm{\ hits \ } \Big(\frac{2}{2-\alpha}\Big)^2p(\infty)^{2-\alpha} \mathrm{\ before \ } 0\Big\} = \left\{ \liminf_{t\to\infty} \xi_t = \infty \right\}, 
\]
so, for $A = (\frac{2}{2-\alpha})^2p(\infty)^{2-\alpha},$
\begin{align}
\label{eq:app.bessel.hits}
\Prob_{\wt\eta_0}\left( \wt\eta_t > 0\right) &< \Prob_{\wt\eta_0}\Big( \wt{\wt\eta}_t > 0 \Big) + \Prob_{\wt\eta_0}\Big( \wt{\wt\eta} \mathrm{\ hits \ } A \mathrm{\ before \ } 0 \Big)  \\  
\notag &= \Prob_{\wt\eta_0}\Big( \wt{\wt\eta}_t > 0 \Big) + \frac{q({\wt\eta}_0) - q(0)}{q(A) -q(0)},
\end{align}
where the scale function $q$ for $\wt{\wt\eta}$ equals, for fixed $\ve_0 >0,$
\[
q(x) = \int_{\ve_0}^{x} \e^{-\frac{\delta}{2} \int_{\ve_0}^y \frac{dz}{z}}dy =  \frac{\ve_0^{\frac{\delta}{2}}}{1-\frac{\delta}{2}} \left( x^{1-\frac{\delta}{2}} -\ve_0^{1-\frac{\delta}{2}}\right).
\]
Since $1-\frac{\delta}{2}=\frac{1}{2-\alpha}$ and 
\begin{align*}
q(A) -q(0) &= \frac{\ve_0^{\frac{\delta}{2}}}{1-\frac{\delta}{2}} A^{1-\frac{\delta}{2}} = \frac{4\ve_0^{\frac{\delta}{2}}}{(1-\frac{\delta}{2})(2-\alpha)^2} p(\wt{C})^{2-\alpha} > 0,
\end{align*}
we get by combining  \eqref{eq:bessel.time} and \eqref{eq:app.bessel.hits} that for some $C=C(t,\alpha,\ve_0, \wt{C})>0,$
\begin{align*}
\Prob_{\wt\eta_0}\left( \wt\eta_t > 0\right) &< C {\wt{\eta}}_0^{\frac{1}{2-\alpha}}, 
\end{align*}
which implies together with \eqref{eq:app.bessel.est} that for some $C_1=C_1(t,\alpha,\wt{C})>0$
\begin{equation}
\label{eq:app.10}
\Prob_x\left( \eta_t > 0 \right) \le C_1 x, \quad x>0, 
\end{equation}
regardless of whether $p(\infty)$ is finite or not.

Define $\wt\theta(a) = \inf\{s\mid \wt\xi_s = a\}, a\ge 0.$ For the diffusion $\xi$ we have
\begin{align*}
\Prob_x\left( \xi_t > 0 \right) &= \Prob_{p(x)}\left( \wt{\xi}_t > 0\right) \\
&\le \Prob_{p(x)} \left( \wt\xi_t > 0, \wt\theta(p(\wt{C})) = \infty\right) + \Prob_{p(x)}\left( \wt\theta(p(\wt{C})) < \wt\theta(0)  \right) \\
&= \Prob_{p(x)} \left( \wt\xi_t > 0, \wt\theta(p(\wt{C})) = \infty\right) + \frac{p(x)}{p(\wt{C})}.
\end{align*}
Since $p(x)\le x, p^{-1}(y)\ge y,$ so  we have for $y\in(0,p(\wt{C}))$
\[ 
\psi(y) \ge  C_\vf \frac{p^{-1}(y)^\alpha}{y^\alpha} \e^{- \int_0^{p^{-1}(y)} \frac{\rho(z)}{1-\vf(z)}dz} \ge C_\vf \e^{- \int_0^{\wt{C}} \frac{\rho(z)}{1-\vf(z)}dz}  =\gamma > 0, 
\]
and  on $\{\theta(p(\wt{C})) = \infty\}$ 
\[
\tau(t) = \int_0^t \frac{ds}{\psi(\wt\xi_{\tau(s)})} \le \gamma^{-1}t. 
\]
Thus we get, using \eqref{eq:app.10},
\begin{align}
\label{eq:add.17}
\Prob_x\left( \xi_t > 0 \right) &\le \Prob_{p(x)} \Big( \wt\xi_t > 0, \tau(s) \le \gamma^{-1}s, s\le \int_0^{\wt\theta} \psi(\wt\xi_r)dr\Big) + \frac{x}{p(\wt{C})} \notag \\
&\le  \Prob_{p(x)} \left( \wt\xi_t = \eta_{\tau^{-1}(t)} > 0, \tau^{-1}(t) \ge \gamma t\right) + \frac{x}{p(\wt{C})} \notag \\
&\le \Prob_{p(x)} \left( \eta_{\gamma t} > 0\right) + \frac{x}{p(\wt{C})} \notag \\
&\le \Big( C_1 + \frac{1}{p(\wt{C})} \Big) x,
\end{align}
which concludes the proof.
\qed

%


\begin{remark}
 \cite[Example 6.4]{KaRuf16Distribution} uses~\cite[Exercise XI.1.22]{ReYor99Continuous} that is stated for only positive dimensions. However, the claim of~\cite[Exercise XI.1.22]{ReYor99Continuous} is known to hold for Bessel processes with arbitrary dimensions (after restricting  to trajectories that do not hit $0$).
\end{remark}

\begin{remark}
In \cite{Mat89Coalescing}, the coalescing property is established and estimates for the number of surviving particles are obtained under weaker assumptions on $\vf$ and for zero drift by studying eigenfunction expansions of the corresponding transitional densities.
\end{remark}

{\it Proof of Lemma \ref{lem:add.1}.} (1) follows from taking expectation in \eqref{eq:add.16}. 

(2) Let $p\ge 2, u_2, u_1\in\mbR, |u_2-u_1|\le 1$ be fixed. 
The flow $X^{\vf,a}$ is a coalescing flow by Theorem~\ref{th:flow.coalescence}. Define
\[
v_t = \sup_{s \in [0,T]} \left| X^{\vf,a}_{0,s}(u_2) - X^{\vf,a}_{0,s}(u_1) \right|, \quad t\in[0,T],
\]
and
\begin{align*}
\theta &= \inf\left\{T; t\in[0,T] \mid  X^{\vf,a}_{0,t}(u_2)= X^{\vf,a}_{0,t}(u_1) \right\}.
\end{align*}
Note that $v_t = v_{\min\{\theta, t\}}, t\in[0,T].$ For some $C_p$
\begin{align*}
v_{t}^p &\le C_p \Big( |u_2-u_1|^p + \int_0^t v_{s}^p ds  +  \sup_{s\in[0,\theta]} \left| m_s \right|^p  \Big),
\end{align*}
where $m$ is a continuous martingale with
\[
\langle m \rangle_t = 2 \int_0^t \left( 1- \vf\left(  X^{\vf,a}_{0,s}(u_2) - X^{\vf,a}_{0,s}(u_1) \right) \right) ds \le 2t, \quad t\in[0,T].
\]
Thus
\begin{align*}
\E v_\theta^p &\le \e^{C_pT} \Big( |u_2-u_1|^p +  \E \sup_{s\in [0,\theta]} \left| m_s \right|^p \Big).
\end{align*}
Let $w$ be a standard Wiener process. For any $a\in[0,T]$ and any $\ve \in (0;1)$
\begin{align*}
 \E \sup_{s\in [0,\theta]} \left| m_s \right|^p  &\le  \E \sup_{s\in [0,\theta]} \left| m_s \right|^p \1\left[ \theta \ge a\right]+  \E \sup_{s\in [0,a]} \left| m_s \right|^p \\
 &\le \Big( \E \sup_{t\in[0,2T]} |w_t|^{\frac{p}{\ve}} \Big)^{\ve} \Prob \left(\theta \ge a \right)^{1-\ve} + \E \sup_{t\in[0,2a]} |w_t|^p.
\end{align*}
Here by \eqref{eq:add.17} and \eqref{eq:bessel.time} for some $C$
\[
\Prob \left(\theta \ge a \right) \le C \frac{|u_2-u_1|}{a^{\frac{1}{2-\alpha}}},
\]
while
\begin{align*}
\E \sup_{t\in[0,2a]} |w_t|^p &= p \int_0^\infty u^{p-1} \Prob\Big( \sup_{t\in[0,2a]} |w_t| \ge u \Big) du \\
&\le \frac{2pa^{1/2} }{\pi^{1/2}} \int_0^\infty u^{p-2} \e^{-\frac{u^2}{4a}} du \\
&= \frac{2p}{\pi^{1/2}} a^{p/2} \int_0^\infty u^{p-2}  \e^{-\frac{u^2}{4}} du.
\end{align*} 
Thus setting $a = |u_2-u_1|^\vk,$ where $\vk  = \gamma(2-\alpha), \gamma \in (0,1),$   gives, for a new constant $C_1$ 
\begin{align*}
\E v_\theta^p &\le C_1 \Big( a^{p/2} + \frac{|u_2-u_1|^{1-\ve}}{a^{\frac{1-\ve}{2-\alpha}}}  \Big) \\
 &= C_1\Big( |u_2-u_1|^{\frac{\gamma (2-\alpha)p}{2}} + |u_2-u_1|^{(1-\ve)(1 - \gamma)}  \Big).
\end{align*}
The optimal $\gamma$ is 
\[
\frac{1-\ve}{1-\ve + \frac{(2-\alpha)p}{2}}.
\]
 \qed

\section{Proofs of Propositions \ref{prop:2.1.1} and \ref{prop:2.1.2}}
 Note that since the flow $X^{\vf,a}$ does not necessarily have the inverse flow one cannot refer to results in~\cite[Chapter 4]{Ku90Stochastic} directly.  

Put
\begin{align*}
\ov{M}_{s,t}(x) &= M_{s,t}(x)-x = X^{\vf, a}_{s,t}(x) - x - \int_s^t a(X^{\vf, a}_{s,r}(x)) dr, \\
\ov{X}^{\vf,a}_{s,t}(x) &= X^{\vf,a}_{s,t}(x) -x, \\
&\quad 0 \le s\le t, x\in\mbR. 
\end{align*}
Each $\ov{M}_{s,\cdot}$ is a standard Wiener process started at $0.$

The following lemma which follows from the definition of a Harris flow is repeatedly used in the sequel. We also skip straightforward calculations.

\begin{lemma}
\label{lem:app.b.1}
For all $s,r,t\ge 0, s \le r \le t,$ and $u,v \in\mbR$
\begin{align*}
\E \ov{M}_{s,t}(u) \ov{M}_{s,t}(v) &= \int_s^t \vf\left(X^{\vf,a}_{s,r}(u) - X^{\vf,a}_{s,r}(v)\right)dr, \\
\E \left( \ov{M}_{s,t}(u) - \ov{M}_{s,t}(v) \right)^2 &= 2\int_{s}^{t} \left( 1 - \vf\left(X^{\vf,a}_{s,r}(u) - X^{\vf,a}_{s,r}(v)\right)\right) dt, \\
\ov{M}_{s,t}(u) &= \ov{M}_{s,r}(x) + \ov{M}_{r,t}\left( X^{\vf,a}_{s,r}(x) \right).
\end{align*}
\end{lemma} 

{\it Proof of Proposition \ref{prop:2.1.1}.} To simplify notation consider $s=0, t=1, t_{n,k} = \frac{k}{2^n}, k=\ov{0,2^n},$
\[
\xi_{n} = \sum_{k=\ov{0,2^n-1}} \ov{M}_{t_{n,k},t_{n,k+1}}(x), \quad n\in\mbN. 
\]
For the sequence $(\xi_n)_{n\in\mbN}$ to be Cauchy in $L_2(\Omega),$ the space of square integrable random variables, it is sufficient to prove 
\begin{equation}
\label{eq:app.b.0}
\alpha_{n,m} = \E \xi_n \xi_m \to t, \quad n,m\to\infty.
\end{equation}
Define 
\begin{align*}
A_{n,k,m} &= \left\{ j\in \{0, \ldots, 2^m-1\} \mid t_{m,j} \in [t_{n,k},t_{n,k+1}] \right\}, \\
&\quad  n,m,k\in\mbN, m\ge n, k=\ov{0,2^n-1}. 
\end{align*}
Assume $m\ge n$ throughout the proof. We have
\begin{align*}
\alpha_{n,m} &= \sum_{k=\ov{0,2^n-1}} \sum_{j\in A_{n,k,m}} \E \ov{M}_{t_{n,k},t_{m,j+1}}(x) \ov{M}_{t_{m,j},t_{m,j+1}}(x) \\
&= \sum_{k=\ov{0,2^n-1}} \sum_{j\in A_{n,k,m}} \E \int_{t_{m,j}}^{t_{m,j+1}} \vf\left( X^{\vf,a}_{t_{m,j},r}(u) - X^{\vf,a}_{t_{m,j},r}(x) \right)\Big\vert_{u=X^{\vf,a}_{t_{n,k},t_{m,j}}(x)} dr \\
&= t - \sum_{k=\ov{0,2^n-1}} \sum_{j\in A_{n,k,m}} \E g_m (x, X^{\vf,a}_{t_{n,k},t_{m,j}}(x)), 
\end{align*}
where 
\begin{align*}
g_m(u,v) =& \E \int_{0}^{\frac{t}{2^m}} \left( 1 - \vf\left(X^{\vf,a}_{0,r}(u) - X^{\vf,a}_{0,r}(v)\right) \right) dt, \\
 &\quad u,v\in\mbR, m\in\mbN.
\end{align*}

Assume $\ve >0$ is fixed. Then
\begin{align*}
g_m(u,v) &\le 2^{-m} \Big[ \Prob\Big( \sup_{r \in[0, 2^{-m}]} \max\left\{ \left| \ov{X}_{0,r}(u)\right|, \left| \ov{X}_{0,r}(v)\right|\right\} \ge \ve \Big) \\
&\qquad\qquad + \sup_{y\in \ov{B}(u-v, 2\ve)} (1-\vf(y)) \Big] \\
&\le 2^{-m} \Big[ 2 \sup_{y\in \{u,v\}} \Prob\Big( \sup_{r\in[0,2^{-m}]} \left| \ov{X}_{0,r}(y)\right| \ge \ve \Big) +  \sup_{y\in \ov{B}(u-v, 2\ve)} (1-\vf(y)) \Big],
\end{align*}
so
\begin{align}
\label{eq:app.b.1}
0 \le t - \alpha_{n,m} &\le 2^{-m+1}  \sum_{k=\ov{0,2^n-1}} \sum_{j\in A_{n,k,m}} \sup_{y\in \left\{x,X^{\vf,a}_{t_{n,k},t_{m,j}}(x)\right\}} \Prob\Big( \sup_{r\in[0,2^{-m}]} \left| \ov{X}_{0,r}(y)\right| \ge \ve \Big) \notag   \\
&\qquad\qquad + 2^{-m+1}  \sum_{k=\ov{0,2^n-1}} \sum_{j\in A_{n,k,m}} \Prob\Big( \big| \ov{X}^{\vf,a}_{t_{n,k},t_{m,j}}(x) \big| \ge \ve \Big) \notag \\
&\qquad\qquad +  \sup_{|y| \le 3\ve} (1-\vf(y)).
\end{align}
Since for some absolute $C>0$ 
\begin{align}
\label{eq:app.b.3}
\sup_{r \in[r_1,r_2]} \left| \ov{X}^{\vf,a}_{r_1, r}(y)\right| &\le C\Big( \sup_{r \in [r_1,r_2]} \left|\ov{M}_{r_1,r}(y)\right| + (1 + |y|)(r_2-r_1)\Big),  \notag \\
&\quad 0 \le r_1 \le r_2 \le 1, y \in \mbR, 
\end{align}
we have, for independent standard Wiener processes $w_1, w_2$ and $j\in A_{n,k,m}$
\begin{multline*}
 \Prob\Big( \left| \ov{X}^{\vf,a}_{t_{n,k},t_{m,j}}(x) \right| \ge \ve \Big) + \sup_{y\in \left\{x,X^{\vf,a}_{t_{n,k},t_{m,j}}(x)\right\}} \Prob\Big( \sup_{r\in[0,2^{-m}]} \left| \ov{X}_{0,r}(y)\right| \ge \ve \Big) \\
  \quad\le 2\Prob\Big(  \sup_{s\in[0, 2^{-n}]} |w_{1s}| + C2^{-m} \sup_{s\in[0, 2^{-n}]} |w_{2s}| \ge \frac{\ve}{C} -2^{-m} -C2^{-n} (1 + |x|)\Big). 
\end{multline*}
Substituting the last estimate into \eqref{eq:app.b.1} implies \eqref{eq:app.b.0}.

It is trivial to check that 
\[
s\mapsto V(s,x) = \left(W_s - W_0, \vf(x-\cdot)\right)_{H_\vf}, 
\] is a martingale w.r.t the filtration of the flow and 
\[
\E V(r_1,x) V(r_2,y) = \min\{r_1, r_2\} \vf(x-y), 
\] 
so it is left to prove that $\{V(t,x) \mid t \ge 0, x\in\mbR\}$ is a jointly Gaussian system. 

Let $u_k = V(r, x_k), x_k\in\mbR, r >0, k=\ov{1,n},  n\in\mbN.$ 
The  Gram--Schmidt orthogonalization produces $v_k \in L_2(\Omega)$ such that
\begin{align*}
\E v_k v_j &= \1\left[k\not= j\right], \\
v_k &= \sum_{i=\ov{1,k}} b_{k,i} u_i,  \\ 
u_k &= \sum_{i=\ov{1,k}} c_{k,i} v_i \quad k,j = \ov{1,n},
\end{align*} 
where the coefficients $\{b_{k,i}\}$ and $\{c_{k,i}\}$ are deterministic and do not depend on $r.$ Define martingales
\[
z_k = \sum_{j=\ov{1,k}} b_{k,j} V(\cdot,x_j), \quad k =\ov{1,n}.
\]
Then 
\[
\langle z_k, z_j \rangle_r = r \E v_k v_j = r \1\left[k\not= j\right], \quad k,j =\ov{1,n}. 
\]
Since uncorrelated continuous martingales with independent increments are independent, the proposition is proved. 
\qed

{\it Proof of Proposition \ref{prop:2.1.2}.} 
By definition of the stochastic integral, 
\[
\int_s^t W\left(X^{\vf,a}_{s,r}(x), dr\right) = L_2-\lim_{n\to\infty} \sum_{k=\ov{0,n-1}} \Big( W(t_{n,k+1},X^{\vf,a}_{s,t_{n,k}}(x)) -W(t_{n,k},X^{\vf,a}_{s,t_{n,k}}(x) \Big),
\]
where $t_{n,k} = s + \frac{k}{n}(t-s), k= \ov{0,n},$ so  
\begin{align}
\label{eq:app.b.30}
I &= \E \Big( \int_s^t W\left(X^{\vf,a}_{s,r}(x), dr\right) - \ov{M}_{s,t}(x) \Big)^2 \notag \\
&= \lim_{n\to\infty}\sum_{k=\ov{0,n-1}} \E \Big[ W\left(t_{n,k+1}, X^{\vf,a}_{s,t_{n,k}}(x)\right) - W\left(t_{n,k}, X^{\vf,a}_{s,t_{n,k}}(x)\right)  \notag \\
& \qquad\qquad\qquad\qquad - \ov{M}_{t_{n,k},t_{n,k+1}}\left(X^{\vf,a}_{s,t_{n,k}}(x)\right) \Big]^2.
\end{align}

Consider uniform partitions of the closed intervals $[t_{n,k},t_{n,k+1}]:$
\[
s^{n,k}_{m,j} = t_{n,k} + \frac{j}{m} (t_{n,k+1} - t_{n,k}), \quad j = \ov{1,m}, k=\ov{0,n-1}.
\]
Then 
\begin{align*}
\ov{M}_{t_{n,k}, t_{n,k+1}}\left( X^{\vf,a}_{s,t_{n,k}}(x)\right) &= \sum_{j=\ov{0,m-1}} \ov{M}_{s^{n,k}_{m,j}, s^{n,k}_{m,j+1}}\Big(X^{\vf,a}_{s,s^{n,k}_{m,j}}(x)\Big), \\
W\left( t_{n,k+1}, X^{\vf,a}_{s,t_{n,k}} \right) - W\left( t_{n,k}, X^{\vf,a}_{s,t_{n,k}} \right) &= L_2\!-\!\lim_{m\to\infty}\! \sum_{j=\ov{0,m-1}} \ov{M}_{s^{n,k}_{m,j}, s^{n,k}_{m,j+1}}\left(X^{\vf,a}_{s,t_{n,k}}(x)\right),
\end{align*}
so returning to \eqref{eq:app.b.30} gives
\begin{align*}
I   &\le \lim_{n\to\infty} \sum_{k=\ov{0,n-1}} \lim_{m\to\infty} \E \Big[ \sum_{j=\ov{0,m-1}} \Big(  \ov{M}_{s^{n,k}_{m,j}, s^{n,k}_{m,j+1}}\Big(X^{\vf,a}_{s,t_{n,k}}(x)\Big)  \\
&\qquad\qquad\qquad\qquad\qquad\qquad - \ov{M}_{s^{n,k}_{m,j}, s^{n,k}_{m,j+1}}\Big(X^{\vf,a}_{s,s^{n,k}_{m,j}}(x)\Big) \Big) \Big]^2 \\
&= \lim_{n\to\infty} \sum_{k=\ov{0,n-1}} \lim_{m\to\infty} \E \sum_{j={0,m-1}} f_{n,m} \Big( X^{\vf,a}_{s,t_{n,k}}(x), X^{\vf,a}_{s,s^{n,k}_{m,j}}(x) \Big),
\end{align*}
where
\begin{align*}
f_{n,m}(u,v) &= 2\E \int_{0}^{\frac{t-s}{nm}} \left( 1 - \vf\left( X^{\vf,a}_{0,r}(u) - X^{\vf,a}_{0,r}(v)\right)\right)dr, \\
&\quad u,v\in\mbR, n,m \in\mbN.
\end{align*}
To obtain the final estimate for \eqref{eq:app.b.30}, one proceeds as in the proof of Proposition~\ref{prop:2.1.1}. 
\qed


\end{appendices}

\section*{Acknowledgments} 
The author is grateful to the anonymous referees for helpful suggestions and corrections which significantly improved the presentation.

This work was supported by a grant from the Simons Foundation (1030291, M.B.V.).
\printbibliography

\end{document}